\documentclass[12pt]{article}
\usepackage{amsmath}
\usepackage{graphicx}
\usepackage{amsfonts}
\usepackage{amssymb}

\newcommand{\gaik}{\gamma_{i,k}}

\newcommand{\Id}{{\rm Id}}

\newcommand{\1}{{\bf 1}}
\newcommand{\la}{\lambda}

\newcommand{\cp}{{\cal P}}

\newcommand{\al}{\alpha}
\newcommand{\ga}{\gamma}

\newcommand{\si}{\sigma}

\newcommand{\E}{{\mathbf E}}

\newcommand{\R}{{\mathbb R}}

\newcommand{\lla}{\left\langle}
\newcommand{\rra}{\right\rangle}
\newcommand{\lcl}{\left\{}
\newcommand{\rcl}{\right\}}
\newcommand{\lp}{\left(}
\newcommand{\rp}{\right)}
\newcommand{\lc}{\left[}
\newcommand{\rc}{\right]}

\newtheorem{theorem}{Theorem}[section]

\newtheorem{corollary}[theorem]{Corollary}

\newtheorem{definition}[theorem]{Definition}

\newtheorem{hypothesis}[theorem]{Hypothesis}
\newtheorem{lemma}[theorem]{Lemma}

\newtheorem{proposition}[theorem]{Proposition}

\newtheorem{remark}[theorem]{Remark}

\begin{document}
\thispagestyle{empty}
 \begin{center}
 {\Large\bf A diluted version of the perceptron model}
\vspace{0.5cm}

\vspace{0.5cm}

{\bf David M\'arquez-Carreras}\footnote{Partially supported by
DGES grant BFM2003-01345.\hfill}, {\bf Carles Rovira}$^1$
and {\bf Samy Tindel}$^2$\\

\vspace{1.5cm}

 $^1$
{\it Facultat de Matem\`atiques,
Universitat de Barcelona,}
\\
\it Gran Via 585,
08007-Barcelona,
Spain
\\
{\it e-mail: davidmarquez@ub.edu, carles.rovira@ub.edu}

\vspace{0.3cm}

$^2$ {\it Institut Elie Cartan, Universit\'e de Nancy 1,}
\\
{\it BP 239,
54506-Vandoeuvre-l\`es-Nancy,
France}
\\
{\it e-mail: tindel@iecn.u-nancy.fr}

 \end{center}

\vspace{2cm}

\begin{abstract}
This note is concerned with a diluted version of the
perceptron model. We establish a replica
symmetric formula at high temperature, which  is achieved
by studying the asymptotic behavior of a given spin magnetization.
Our main task
will be to identify the order parameter of the system.

\end{abstract}

\vspace{2cm}


\vspace{0.3cm}


\newpage

\section{Introduction}

A wide number of spectacular advances have occurred in the spin glasses
theory during the last past years, and it could easily be argued that this
topic, at least as far as the Sherrington-Kirkpatrick model is concerned,
has reached a certain level of maturity from the mathematical point of view:
the cavity method  has been set in a clear and effective way in \cite{Tb},
 some monotonicity properties  along a smart path have been discovered
 in \cite{GT}, and these elements have been combined in \cite{Tam} in order to
 obtain a completely rigorous proof of the Parisi solution \cite{MPV}.

\vspace{0.3cm}

However, there are some canonical models of mean field spin glasses for
which the basic theory is far from being complete, and this paper proposes
 to study the high temperature behavior of one of them, namely the diluted
 perceptron model, which can be described as follows: for $N\ge 1$, consider
 the configuration space $\Sigma_N=\{-1,1\}^N$, and for
$\mathop{\pmb{\si}}=(\si_1,\dots,\si_N)\in \Sigma_N$, define a Hamiltonian
$-H_{N,M}(\mathop{\pmb{\si}})$ by
\begin{equation}\label{ehamin}
-H_{N,M}
(\mathop{\pmb{\si}})= \sum_{k \le M} \eta_k \ u \left( \sum_{i \le
N} g_{i,k}\ \ga_{i,k}\ \si_{i} \right).
\end{equation}
In this Hamiltonian, $M$ stands for a positive integer such that
$M=\alpha N$ for a given $\alpha\in (0,1)$; $u$ is a bounded
continuous function defined on $\mathbb{R}$; $\{g_{i,k},\ i\ge 1,
k\ge 1\}$ and $\{\gaik,\  i\ge 1, k\ge 1\}$ are two independent
families of independent random variables, $g_{i,k}$ following  a
standard Gaussian law and $\gaik$ being a Bernoulli random
variable with parameter $\frac{\ga}{N}$, which we denote by
$B(\frac{\ga}{N})$. Eventually, $\{\eta_{k},\ k\ge 1\}$ stands for
an arbitrary family of numbers, with $\eta_k\in\{0,1  \}$, even
if the case of interest for us will be $\eta_k=1$ for all $k\le M$.
Associated to this Hamiltonian, define a random Gibbs measure
$G_N$ on $\Sigma_N$, whose density with respect to the uniform
measure $\mu_N$ is given by $Z_{N,M}^{-1}
\exp\left(-H_{N,M}(\mathop{\pmb{\si}})\right)$, where the
partition function $Z_{N,M}$ is defined by
$$
Z_{N,M}=\sum_{\mathop{\pmb{\si}}\in \Sigma_N}
\exp\left(-H_{N,M}(\mathop{\pmb{\si}})\right).
$$
In the sequel, we will denote by $\langle f \rangle$ the average of a function
$f:\Sigma_N^n\longrightarrow \mathbb{R}$ with respect to
$dG_N^{\otimes n}$, i.e.
$$
\langle
f\rangle=Z_{N,M}^{-n}\sum_{(\mathop{\pmb{\sigma}}^1,\ldots,\mathop{\pmb{\sigma}}^n)\in\Sigma_N^n}
f(\mathop{\pmb{\sigma}}\!{^1},\ldots,\mathop{\pmb{\sigma}}\!{^n})\exp\lp
-\sum_{l\le n}H_{N,M}(\mathop{\pmb{\sigma}}\!{^l}) \rp.
$$


\vspace{0.3cm}

The measure described above is of course a generalization of the usual perceptron model,
which has been introduced for neural computation purposes
(see \cite{HKP}), and whose high temperature behavior has been described in
\cite[Chapter 3]{Tb}, or \cite{ST} for an approach based on convexity
properties of the Hamiltonian. 
Indeed the usual perceptron model is induced by a Hamiltonian
$\hat H_{N,M}$ on $\Sigma_N$ given by
\begin{equation}\label{ndper}
-\hat H_{N,M}
(\mathop{\pmb{\si}})= \sum_{k \le M}
u \left( \frac{1}{N^{1/2}}\sum_{i \le
N} g_{i,k}\  \si_{i} \right),
\end{equation}
where we have kept the notations introduced for equation (\ref{ehamin}).
Thus, our model can be seen as a real diluted version of (\ref{ndper}), in the  sense
 that in our model, each condition $\sum_{i\le N}g_{i,k}\ga_{i,k}\si_{i}\ge 0$ only involves, in average,
  a finite number of spins, uniformly in $N$. It is worth noticing at
that
 point that this last requirement fits better to the initial neural
computation
  motivation, since in a one-layer perceptron, an output is generally
obtained by
  a threshold function applied to a certain number of spins, that does
not grow
  linearly with the size of the system.
Furthermore, our coefficient $\gamma$ is arbitrarily large, which
means that the global interaction between spins is not trivial.
Another motivation for the study
of the system
  induced by (\ref{ehamin}) can be found in \cite{CL}. Indeed, in this
latter article,
  a social interaction model is proposed, based on a Hopfield-like
(or perceptron-like)
  diluted Hamiltonian with parameter $N$ and $M$, where $N$ represents
the number of
  social agents, and $M$ the diversity of these agents, the number of
interactions
  of each agent varying with the dilution parameter. However, in \cite{CL},
the
  equilibrium of the system is studied only when $M$ is a fixed number.
The result
   we will explain later on can thus be read as follows: as soon as the
diversity
   $M$ does not grow faster than a small proportion of $N$,
the capacity
the social interaction system is not attained

\vspace{0.3cm}

Let us turn now to a brief description of the results contained in this paper:
in fact, we will try to get a replica symmetric formula for the system when
$M$ is a small proportion of $N$, which amounts to identify  the limit of
$\frac1N\log(Z_{N,M})$ when $N\to\infty$, $M=\al N$. This will be achieved, as
in the diluted SK model studied through the cavity method (see \cite{FT} for a
study based on monotonicity methods), once the limiting law for the magnetization
$\langle \si_i\rangle$ is obtained. This will thus be our first aim, and
in order to obtain that  result, we will try to adapt the method designed in
\cite[Chapter 7]{Tb}. However, in our case, the identification of the
limiting law for $\langle \si_i\rangle$ will be done through an intricate fixed
point argument, involving a map $T:{\bf P}\to{\bf P}$ (where ${\bf P}$ stands
for the set of probability measures on $[-1,1]$), which in turn  involves a kind of
$\cp(\la)^{\otimes\cp(\mu)}$ measure, for two independent Poisson measures
$\cp(\la)$ and $\cp(\mu)$. For sake of readability, we will give the details of
(almost) all the computations we will need in order to establish our replica
symmetric formula, but it should be mentioned at that point that our
main contribution, with respect to \cite[Chapter 7]{Tb}, is that construction
of the invariant measure.

\vspace{0.3cm}

More specifically, our paper is divided as follows:
\begin{itemize}
\item
At Section \ref{spincor}, we will establish a decorrelation result for two
arbitrary spins. Namely, setting $U_{\infty}=\|u\|_{\infty}$, for
$\al U_{\infty}$ small enough, we will show that
$$
\E\lc\left|\langle \si_1 \si_2\rangle - \langle \si_1\rangle \langle
\si_2\rangle\right|\rc \le \frac{K}{N},
$$
for a constant $K>0$.
\item
At Section \ref{magn}, we will study the asymptotic behavior of the
magnetization of $m$ spins, where $m$ is an arbitrary integer.
Here again, if $\al U_{\infty}$ is small enough, we will see that
$$
\E \lc \sum_{i\le m} \left\vert \langle
\si_i\rangle-z_i\right\vert \rc \le \frac{K m^3}{N},
$$
where $z_1,\dots,$ $z_m$ is a family of i.i.d random variable, with law
$\mu_{\al,\gamma}$, and $\mu_{\al,\gamma}$ is the fixed point of the map
$T$ alluded to above, whose precise description will be given at
the beginning of Section \ref{magn}.
\item
Eventually, at Section \ref{symf}, we obtain the replica symmetric formula for our model: set
\begin{eqnarray*}
\bar V_p &=&
\int \Big \langle \exp( u(\sum_{i \le p} g_{i,M}
\sigma_i )) \Big \rangle_{(x_1,\ldots,x_{p})}
d\mu_{\alpha,\gamma}(x_1)  \cdots
d\mu_{\alpha,\gamma}(x_{p})\\
G(\gamma) &=&
\alpha \log \left( \sum_{p=0}^\infty \exp(-\gamma)
\frac{\gamma^p}{p!} \E \left[ \frac{\bar V_{p+1}}{\bar V_p}
\right] \right),
\end{eqnarray*}
where $\langle \cdot \rangle_x$  means integration with respect to
the product measure $\nu$ on $\{-1,1\}^p$ such that
$\int \sigma_i d\nu = x_i$. Let $F:[0,1]\to\R^+$ be defined by
$F(0)=\log 2 - \alpha u(0)$ and $F'(\gamma)=G(\gamma)$. Then,
if $\al U_{\infty}$ is small enough, we will get that
$$
\left\vert \frac{1}{N}\E\lc \log(Z_{N,M}) \rc
- F (\gamma) \right\vert \le \frac{K}{N},
$$
for a strictly positive constant $K$.
\end{itemize}
All these results will be described in greater detail in the
corresponding sections.

\section{Spin correlations}\label{spincor}

As in \cite[Chapter 7]{Tb}, the first step towards a replica
symmetric formula will be to establish a decorrelation result for
two arbitrary spins in the system. However, a much more general
property holds true, and we will turn now to its description: for
$j\le N$, let $T_j$ be the transformation of $\Sigma_N^n$ that,
for a configuration
$(\mathop{\pmb{\si}}^1,\dots,\mathop{\pmb{\si}}^n)$ in
$\Sigma_N^n$, exchanges the $j$-th coordinates of
$\mathop{\pmb{\si}}^1$ and $\mathop{\pmb{\si}}^2$. More
specifically, let $f:\Sigma_N^n\to\R$, with $n\geq2$, and let us
write, for $j\leq N$,
$$
f=f \left(  \mathop{\pmb{\si}_{j^c}^{1}},\sigma_{j}^{1};
\mathop{\pmb{\si}_{j^c}^{2}},\sigma_{j}^{2};\ldots;
\mathop{\pmb{\si}_{j^c}^{n}} ,\sigma_{j}^{n}\right) ,
$$
where, for $l=1,\dots,n$,
$\mathop{\pmb{\si}_{j^c}^{l}}=(\si_1^l\dots,\si_{j-1}^l,\si_{j+1}^l,\dots,
\si_N^l)$. Then define $f\circ T_j$ by
\begin{equation}\label{defuj}
f\circ T_j(\mathop{\pmb{\si}^{1}},\ldots,\mathop{\pmb{\si}^{n}})
=
f
\left(  \mathop{\pmb{\si}_{j^c}^{1}},\sigma_{j}^{2};
\mathop{\pmb{\si}_{j^c}^{2}},\sigma_{j}^{1};\ldots;
\mathop{\pmb{\si}_{j^c}^{n}}
,\sigma_{j}^{n}\right).
\end{equation}
For $j\le N-1$,
we will call $U_j$ the equivalent
transformation on $\Sigma_{N-1}^n$.
\begin{definition}\label{dproperty}
We say that Property \textbf{P}$(N,\ga_0,B)$ is satisfied
if the following requirement is true: let $f$ and $f'$ be two
functions on $\Sigma_N^n$ depending on $m$ coordinates, such that
$f \ge 0$, $f' \circ T_N = -f'$,  and there exists $Q \ge 0$
such that $\vert f'\vert \le Q f$; then if  $\ga \le \ga_0$ we
have
$$
\E \left\vert \frac{\langle f' \rangle}{\langle f \rangle}
\right\vert \le \frac{m Q B}{N},
$$
for any Hamiltonian of the form (\ref{ehamin}), uniformly in $\eta$.
\end{definition}
Set now $U_\infty=\|u\|_\infty$. With Definition
\ref{dproperty} in hand, one of the purposes of this section
is to prove the following Theorem.

\begin{theorem}\label{tproperty}
Let $\gamma_0$ be a positive number, and $U_\infty$ be small enough, so that
\begin{equation}\label{e712}
 4 U_\infty \ \al \ga_0^2\ e^{4U_\infty}\ e^{\al \ga_0
(e^{4U_\infty}-1)}\ \left(3+ 2 \ga_0 + \al (\ga_0^2 + \ga_0^3)
e^{4U_\infty} \right) <1.
\end{equation}
Then there exists a number
$B_0(\ga_0,U_\infty)$ such that if $\ga\le \ga_0$,
the property
\textbf{P}$(N,\ga_0,B_0)$ holds true for each $N\ge 1$.
\end{theorem}
In the previous theorem, notice that the value of $\gamma_0$ has been picked arbitrarily. Then we have to choose $U_\infty$, which also contains implicitly the temperature parameter, accordingly. Let us also mention that the spin decorrelation follows easily from the last result:
\begin{corollary}
Assuming (\ref{e712}) there exists $K>0$ such that,
for all $\gamma<\gamma_0$,
$$\E\left|\langle \si_1 \si_2\rangle - \langle \si_1\rangle \langle
\si_2\rangle\right| \le \frac{K}{N}.$$
\end{corollary}

\vspace{0.5cm}

\noindent {\bf Proof:} It is an easy consequence of property
\textbf{P}$(N,\ga_0,B_0)$ applied to $n=2$, $f=1$ and
$f'(\mathop{\pmb{\sigma}}^1,\mathop{\pmb{\sigma}}^2)= \sigma_1^1(
\si_2^1-\si_2^2)$.

\hfill $\square$

\medskip

We will prepare now the ground for the proof of Theorem \ref{tproperty},
which will be based on  an induction argument
 over $N$.
A first step in this direction will be to state the cavity formula
for our model: for $\mathop{\pmb{\sigma}}=(\si_1,\dots,\si_N)\in
\Sigma_N$, we write
$$
\mathop{\pmb{\rho}}\equiv(\rho_1,\dots,\rho_{N-1})
= (\si_1,\dots,\si_{N-1})\in\Sigma_{N-1}.
$$
Then the Hamiltonian (\ref{ehamin}) can be decomposed into
\begin{equation*}
-H_{N,M} (\mathop{\pmb{\si}}) = \sum_{k \le M} \eta_k \ga_{N,k} u
\left( \sum_{i \le N-1} g_{i,k} \ga_{i,k} \si_{i} + g_{N,k} \si_N
\right) - H_{N-1,M}^-(\mathop{\pmb{\rho}}),
\end{equation*} with
\begin{equation}\label{enoha} -H_{N-1,M}^-(\mathop{\pmb{\rho}})= \sum_{k \le M} \eta_k^- u
\left( \sum_{i \le N-1} g_{i,k} \ga_{i,k} \si_{i}
\right),
\quad\mbox{ and }\quad
\eta_k^-=\eta_k(1-\gamma_{N,k})
.
\end{equation}
Note that in $H_{N-1,M}^-$, the coefficients
$\eta_k^-=\eta_k (1-\ga_{N,k})$ are
not deterministic, and hence $H_{N-1,M}^-$ is not really of
the same kind as $H_{N,M}$.
However, this problem can be solved by conditioning on
$\{\ga_{N,k}, k\le M\}$. Then, given the randomness contained in
the $\ga_{N,k}$, the expression $H_{N-1,M}^-(\mathop{\pmb{\rho}})$
is a Hamiltonian of a $(N-1)$-spin system with $\ga_{i,k} \sim
B(\frac{\ga^-}{N-1})$, where $\ga^-=\ga \frac{N-1}{N}$ and so
$\ga^-\le \ga \le \ga_0$.

Thus, given  a function $f: \Sigma_N^n \longrightarrow \R,$ we
easily get the following decomposition of the mean value of $f$
with respect to $G_N^{\otimes n}$:
\begin{equation}\label{dcpm}
\langle f \rangle = \frac {\langle {\bf Av} f \xi
\rangle_-}{\langle {\bf Av} \xi \rangle_-},
\end{equation}
with
\begin{equation}\label{defxi}
\xi= \exp \left( \sum_{l \le n} \sum_{k \le M} \eta_k \ga_{N,k} u
\left( \sum_{i \le N-1} g_{i,k} \ga_{i,k} \si_{i}^l
 + g_{N,k} \si_N^l \right) \right),
\end{equation}
and with $\langle \bar f \rangle_-$ defined, for a given $\bar
f:\Sigma_{N-1}^n\to\R$, by
 $$
 \langle \bar f \rangle_-
 =
 \frac{ \sum_{(\mathop{\pmb{\rho}}^1,\ldots,\mathop{\pmb{\rho}}^n) \in \Sigma_{N-1}^n}
\bar f (\mathop{\pmb{\rho}}^1,\ldots,\mathop{\pmb{\rho}}^n)
\exp\left( - \sum_{l \le n} H_{N-1,M}^-
(\mathop{\pmb{\rho}}^l)\right)}{
\sum_{(\mathop{\pmb{\rho}}^1,\ldots,\mathop{\pmb{\rho}}^n) \in
\Sigma_{N-1}^n} \exp\left( - \sum_{l \le n} H_{N-1,M}^-
(\mathop{\pmb{\rho}}^l)\right)}.
$$
Notice also that in expression (\ref{dcpm}), ${\bf Av}$ stands for
the average with respect to the last component of the system,
namely if
$f=f(\mathop{\pmb{\rho}}^1,\si_N^1,\ldots,\mathop{\pmb{\rho}}^n,\si_N^n)$,
then
$$
{\bf Av} f(\mathop{\pmb{\rho}^1}\ldots,\mathop{\pmb{\rho}^n}) =
\frac{1}{2^n}\sum_{\si_N^j=\pm 1,j\le n}
f(\mathop{\pmb{\rho}^1},\si_N^1,\ldots,\mathop{\pmb{\rho}^n},\si_N^n).
$$

\vspace{0.3cm}

Let us introduce now a little more notation:
in the sequel we will have to take expectations  for a fixed value of  $\xi$
given at (\ref{defxi}). Let us denote thus
by $\E_{\ga_N}$ the expectation given $\ga_{N,k},\ k \le M$, and
define
\begin{equation}\label{emgn}
\E_{-,\ga_N}[\ \cdot\ ] = \E_{\ga_N} \left[\ \cdot\  \vert \
g_{N,k},\ g_{i,k},\ \ga_{i,k},\ i \le N-1, \ k \in D^M_{N,1}\
\right],
\end{equation}
where $D^M_{N,1}$ is given by
$$
D^M_{N,1} = \{ k \le M; \, \ga_{N,k} = 1 \}.
$$
One has to be careful about the way all these conditioning are
performed, but it is worth observing that the set $D^M_{N,1}$ is
not too large: indeed, it is obvious that, setting $|A|$ for the
size of a set $A$, we have
\begin{equation}\label{defr1}
\vert D^M_{N,1}  \vert = \sum_{k \le M } \ga_{N,k},
\end{equation}
and thus
$$
\E
\vert D^M_{N,1}  \vert = M \frac{\ga}{N}=\al \ga.
$$

\vspace{0.3cm}

Let us go on now with the first step of the induction procedure for the
proof of Theorem \ref{tproperty}:
in \textbf{P}$(N,\ga_0,B)$ we can assume without loss of
generality that $f$ and  $f'$ depend on the coordinates
$1,\ldots,m-1,N.$ Moreover, since $\vert f' \xi \vert \le Q f
\xi$, we have
$$ \vert \langle {\bf Av} f' \xi \rangle_- \vert \le
\langle {\bf Av} \vert f' \xi \vert \rangle_{-} \le \langle Q {\bf
Av} f \xi \rangle_-,$$ and hence
\begin{equation}\label{efitaA}
\left\vert \frac{\langle {\bf Av} f' \xi \rangle_-}{\langle {\bf
Av} f \xi \rangle_-} \right\vert \le Q.\end{equation}
We now define the following two events:
\begin{eqnarray*}
\Omega_1 & = & \{ \exists p \le m-1, k \in D^M_{N,1};\,\ga_{p,k}=1
\}
\\
& = & \{ \exists p \le m-1, k \le M;\,\ga_{p,k}=\ga_{N,k}=1 \},
\\[2mm]
\Omega_2 & = & \{ \exists j  \le N-1, k_1, k_2 \in D^M_{N,1};\,
\ga_{j,k_1}=\ga_{j,k_2}=1 \}
\\
& = & \{ \exists j  \le N-1, k_1, k_2 \le M;\,
\ga_{j,k_1}=\ga_{j,k_2}=\ga_{N,k_1}=\ga_{N,k_2}=1 \}.
\end{eqnarray*}
These two events can be considered as exceptional. Indeed, it is readily checked that
$$ P(\Omega_1) \le \al
\frac{\ga^2}{N} (m-1), \qquad P(\Omega_2) \le \al^2 \ga^4
\frac{N-1}{N^2}.$$
Thus, if $\Omega =  \Omega_1 \cup \Omega_2$, we get
$$
P(\Omega)\le \frac{\al \ga^2 (m-1) + \al^2 \ga^4}{N},
$$
and using this fact together with (\ref{efitaA}), we have
\begin{eqnarray}
\E \left\vert \frac{\langle f' \rangle}{\langle f \rangle}
\right\vert & = & \E \left\vert \frac{\langle {\bf Av} f' \xi
\rangle_-}{\langle {\bf Av} f \xi \rangle_-} \right\vert\nonumber \\
& = & \E \left( \1_{\Omega} \left\vert \frac{\langle {\bf Av} f'
\xi \rangle_-}{\langle {\bf Av} f \xi \rangle_-} \right\vert
\right) + \E \left( \1_{\Omega^c} \left\vert \frac{\langle {\bf
Av} f' \xi \rangle_-}{\langle {\bf Av} f \xi \rangle_-}
\right\vert \right)\nonumber \\
& \le & Q \frac{\al \ga^2 (m-1)+ \al^2 \ga^4}{N} + \E \left(
\1_{\Omega^c} \left\vert \frac{\langle {\bf Av} f' \xi
\rangle_-}{\langle {\bf Av} f \xi \rangle_-} \right\vert
\right).\label{edespr}
\end{eqnarray}
Consequently, in order to prove Theorem \ref{tproperty} we only need to bound
accurately the expectation of the right-hand side of
(\ref{edespr}) by means of  the induction hypothesis.
To this purpose, we will introduce some new notations and go through a series
of lemmas: set
$$
J_1= \{ j\le N;\, \ga_{j,k}=1 \, {\rm for \, some \,} k \in D^M_{N,1}
\}- \{N \},
$$
and observe that, when $\Omega_1$ does not occur,
$$J_1 \cap
\{1,\ldots,m-1\}=\emptyset.
$$
Denote $\vert J_1 \vert ={\rm card}
(J_1)$ and write an enumeration of $J_1$ as follows:
$J_1=\{j_1,\ldots,j_{\vert J_1 \vert} \}.$
\begin{lemma}\label{l724}
Let $U_j$ be the transformation  defined at (\ref{defuj}),
and $f':\Sigma_N^n\to\R$ such that $f'\circ T_N=-f'$.
When  $\Omega_1$ does not occur, we have
$$
({\bf Av} f' \xi ) \circ \prod_{j \in J_1} U_j  = - {\bf Av} f'
\xi.$$
\end{lemma}

\vspace{0.5cm}

\noindent {\bf Proof:} The proof of this lemma can be done
following the steps of \cite[Lemma 7.2.4]{Tb}, and we include it
here for sake of readability. Set $T= \prod_{j \in J_1} T_j.$
Since $f'$ depends only on the coordinates $\{1,\ldots,m-1,N\}$
and this set is disjoint from $J_1$, we have  $f' \circ T = f'$.
Moreover,
$$f' \circ T \circ T_N = f' \circ T_N = -f'.$$
On the other hand, $\xi$ only depends on $J_1 \cup \{N\}$ and
using
$$
\xi(\si^1,\si^2,\ldots,\si^n)=\xi(\si^2,\si^1,\ldots,\si^n),
$$
we obtain
$$
\xi \circ T \circ T_N = \xi.$$
Hence
$$
(f' \xi ) \circ T \circ T_N = -f' \xi,
$$
and, since $T_N^2=\Id$, we get
\begin{equation}\label{asim1}
(f'\xi) \circ T = -(f' \xi) \circ T_N.
\end{equation}
Eventually,
\begin{eqnarray}
{\bf Av} [ (f'\xi) \circ T_N ] & = & {\bf Av} f' \xi, \label{sim1}\\
{\bf Av} [ (f'\xi) \circ T ] & = & ({\bf Av} f' \xi) \circ
\prod_{j \in J_1} U_j.\label{sim2}
\end{eqnarray}
The proof is now easily concluded by plugging (\ref{sim1})
and (\ref{sim2}) into (\ref{asim1}).

\hfill $\square$

\medskip

Let us now go on with the proof of Theorem \ref{tproperty}: thanks to Lemma
\ref{l724}, when $\Omega_1$ does not occur, we can write
\begin{equation}\label{edefifs}
{\bf Av} f' \xi= \frac12 \left[ {\bf Av} f' \xi - ({\bf Av} f'
\xi) \circ \prod_{s \le \vert J_1 \vert} U_{j_s} \right] = \frac12
\sum_{1 \le s \le \vert J_1 \vert} f_s,\end{equation} with
$$
f_s =({\bf Av} f' \xi) \circ \prod_{l \le s-1} U_{j_l} - ({\bf Av}
f' \xi) \circ \prod_{l \le s} U_{j_l}.
$$
Notice that $U_j^2= \Id$, and that $f_s$ enjoys the same kind of
antisymmetric property as $f'$, since $f_s \circ U_{j_s} = - f_s.$

\vspace{0.3cm}

Define $R_1 = \vert D^M_{N,1} \vert$. Then, recalling relation (\ref{defr1}),
we have
$$
R_1  = \vert D^M_{N,1} \vert= \sum_{k \le M} \ga_{N,k},
$$
and let us enumerate as $k_1,\ldots,k_{R_1}$ the values $k\le M$ such that $\ga_{N,k}=1.$
We also define $I_1^1,\ldots,I^1_{R_1}$ as follows:
$$
I^1_v = \{ j\le N-1;\, \ga_{j,k_v}=1 \}, \quad\mbox{ for }\quad v\le
R_1,
$$ and observe that we trivially have
\begin{equation}\label{easterisc}
J_1=\bigcup_{v\le R_1} I^1_v.\end{equation}
Moreover, when
$\Omega_2$ does not occur, we have
$$I^1_{v_1} \cap I^1_{v_2} = \emptyset, \quad {\rm if} \ v_1\neq
v_2.$$
Then, on $\Omega^c$, we get
\begin{equation}\label{ecardJ}
|J_1|={\rm Card} (J_1)= \sum_{v\le R_1}|I^1_v|.
\end{equation}
Furthermore, it is easily checked that, for each $v$,
and conditionally on the $\ga_{N,k}$,
the quantity
$|I^1_v|$ is a binomial random variable with parameters $N-1$ and
$\frac{\gamma}{N}$, which we denote by
${\rm Bin}( N-1,\frac{\gamma}{N})$.

\vspace{0.3cm}

With all these notations in mind, our next step will be to
 bound $f_s$ in function of $f$, in order to get a similar condition to that of Definition
\ref{dproperty}:
\begin{lemma}\label{l725}
Recall that $U_\infty=\| u \|_\infty$.
Then, on  $\Omega^c$, for $j_s \in I^1_v$, we have
$$
\vert f_s \vert \le  \hat Q {\bf Av} f \xi,
$$
where
$$
\hat Q \equiv 4 Q U_\infty \exp \left( 4 U_\infty R_1 \right).
$$
\end{lemma}

\vspace{0.5cm}

\noindent {\bf Proof:} Let us decompose $\xi$ as $\xi = \xi'\ \xi''$, with
\begin{eqnarray}
\xi' & =  & \exp \left( \sum_{3 \le l \le n} \sum_{k \le M} \eta_k
\ga_{N,k} u \left( \sum_{i \le N-1} g_{i,k} \ga_{i,k} \si_{i}^l
 + g_{N,k} \si_N^l \right) \right),\nonumber\\
\xi'' & =  & \exp \left( \sum_{l \le 2} \sum_{k \le M} \eta_k
\ga_{N,k} u \left( \sum_{i \le N-1} g_{i,k} \ga_{i,k} \si_{i}^l
 + g_{N,k} \si_N^l \right) \right).\label{defxis}
 \end{eqnarray}
Thus
 \begin{eqnarray*}
 \xi & \ge & \xi' \exp  \left( - \sum_{l \le 2} \sum_{\bar v \le R_1} \left\vert u \left( \sum_{i \le N-1}
  g_{i,k_{\bar v}} \ga_{i,k_{\bar v}} \si_{i}^l
 + g_{N,k_{\bar v}} \si_N^l \right) \right\vert \right)\\
 & \ge & \xi' \exp \left( -2 U_\infty R_1 \right),
 \end{eqnarray*}
 and hence
 \begin{equation}\label{e730}
 {\bf Av} f \xi \ge ({\bf Av} f \xi' )
\exp \left( -2 U_\infty R_1  \right).
\end{equation}
On the other hand, since $f'$ only depends on
$\{1,\ldots,m-1,N\}$, we have $f' \circ T_{j_l}=f'$ for any
$l\le |J_1|$ on $\Omega^c$, which yields
\begin{eqnarray}
 f_s & = &
({\bf Av} f' \xi) \circ \prod_{l \le s-1} U_{j_l} - ({\bf Av} f'
\xi) \circ \prod_{l \le s} U_{j_l} \nonumber \\
& = & {\bf Av} \left( (f' \xi) \circ \prod_{l \le s-1} T_{j_l} -
(f' \xi) \circ \prod_{l \le s} T_{j_l} \right)  \nonumber \\ & = &
{\bf Av} \left( f' \left( \xi \circ \prod_{l \le s-1} T_{j_l} -
\xi \circ \prod_{l \le s} T_{j_l} \right)\right), \label{44}
\end{eqnarray}
where we have used the fact that $J_1$ can be written as
$J_1=\{j_1,\ldots,j_{\vert J_1 \vert} \}$. Moreover, for
any $l$, by construction of $\xi'$,  we have $\xi' \circ T_{j_l} = \xi'.$
Thus,
\begin{equation}\label{441} \xi \circ \prod_{l \le s-1} T_{j_l} -  \xi \circ \prod_{l
\le s} T_{j_l} = \xi' \left[\xi'' \circ \prod_{l \le s-1} T_{j_l}
-  \xi'' \circ \prod_{l \le s} T_{j_l} \right].
\end{equation}
Set now
$$
\Gamma=\sup_\si \left\vert \xi'' \circ \prod_{l \le s-1} T_{j_l} -
\xi'' \circ \prod_{l \le s} T_{j_l} \right\vert  = \sup_\si \vert
\xi'' - \xi'' \circ T_{j_s} \vert.
$$
Then, from (\ref{44}) and (\ref{441}), and invoking the fact that
$\vert f' \vert \le Q f$, we get
\begin{equation}\label{33}
\vert f_s \vert \le \Gamma {\bf Av} ( \vert f' \vert \xi' ) \le Q
\Gamma {\bf Av} f \xi'.
\end{equation}
We now bound $\Gamma$: recall that $\xi''$ is defined by (\ref{defxis}), and thus
$$
\xi''=\prod_{\bar v \le R_1} \xi_{\bar v},
$$
with
$$\xi_{\bar v} = \exp \left( \sum_{l \le 2} \eta_{k_{\bar v}}  u \left(  \sum_{i \le N-1} g_{i,k_{\bar v}} \ga_{i,k_{\bar v}}
\si_{i}^l
 + g_{N,k_{\bar v}} \si_N^l \right) \right).
$$
Recall now that we have assumed that $j_s\in I_v^1$. Therefore, we have
$j_s \notin I^1_{\bar v}$ if $\bar v\neq v$, according to the fact that
$I^1_{v} \cap I^1_{\bar v} = \emptyset$ on $\Omega^c$. Hence
$$
 \xi_{\bar v} \circ T_{j_s} = \xi_{\bar v},
$$
and
\begin{equation}\label{xisot}
 \xi'' - \xi'' \circ T_{j_s} = ( \xi_v - \xi_v \circ T_{j_s} )
 \prod_{\bar v \not= v} \xi_{\bar v}.
\end{equation}
On the other hand, since $\vert e^x - e^y \vert \le \vert x-y \vert e^a$
 for  $\vert x \vert, \vert y \vert \le a$, we obtain
\begin{equation}\label{difxiv}
 \vert \xi _v - \xi_v \circ T_{j_s} \vert \le
4 U_\infty e^{2 U_\infty },
\end{equation}
and we also have the trivial bound
\begin{equation}\label{trivxi}
\xi_{\bar v} \le e^{2 U_\infty }.
\end{equation}
Thus, plugging (\ref{difxiv}) and (\ref{trivxi}) into (\ref{xisot}), we get
 $$
\Gamma \le 4 U_\infty e^{2 U_\infty R_1}.
$$
Combining this bound with (\ref{e730}) and (\ref{33}),
the proof is now easily completed.

\hfill$\square$

\bigskip

We are now ready to start the induction procedure  on
$\textbf{P}(N,\ga_0,B)$, which will use the following elementary lemma
(whose proof
 is left to the reader).
\begin{lemma}\label{lbin} Let $R$ be a random variable following the
${\rm Bin}(M,\frac{\ga}{N})$ distribution, and $\la$ be a positive number.
Then
\begin{eqnarray}
\E\left[ R e^{\la R}\right]&\le&  \al \ga e^\la e^{ \al \ga  (e^\la -1)},\label{ebin1}\\
\E\left[ R^2 e^{\la R}\right] &\le& \al^2  \ga^2 e^{2\la} e^{ \al
\ga (e^\la -1)} + \al \ga e^\la e^{ \al \ga  (e^\la
-1)}.\label{ebin2}
\end{eqnarray}
\end{lemma}
Let us proceed now with the main step of the induction:
\begin{proposition}\label{p726}
Assume that \textbf{P} $(N-1,\ga_0,B)$ holds for $N\ge 2$ and $\ga
\le \ga_0$. Consider $f$ and $f'$ as in Definition
\ref{dproperty}. Then
\begin{equation}\label{e733}
\E\lc\left\vert
\frac{\langle f' \rangle}{\langle f \rangle} \right\vert\rc
\le
\frac{m Q }{N}\left(\al \ga^2 + \al^2 \ga^4+4B\Upsilon(\alpha,
\ga, U_\infty)\right),
\end{equation} where
\begin{displaymath}\Upsilon(\alpha, \ga, U_\infty)=U_\infty \al \ga^2\
e^{4U_\infty}\ e^{\al \ga (e^{4U_\infty}-1)}\ \left(3+ 2\ga + \al
(\ga^2 + \ga^3) e^{4U_\infty} \right).
\end{displaymath}
\end{proposition}

\vspace{0.5cm}

\noindent {\bf Proof:}  Using (\ref{edespr}) and (\ref{edefifs}),
we have
$$
\E \lc\left\vert \frac{\langle f' \rangle}{\langle f \rangle}
\right\vert\rc
\le  Q \frac{\al \ga^2 (m-1) + \al^2 \ga^4}{N} +
\frac12 \E \left( \1_{\Omega^c} \sum_{s \le \vert J_1 \vert}  \frac{
\vert \langle f_s \rangle_- \vert}{\langle {\bf Av} f \xi
\rangle_-} \right).
$$
However, on $\Omega^c$, the functions $f_s$  and  ${\bf Av} f \xi$ depend
on $m-1+|J_1|$ coordinates. Since $\ga^-\le \ga$ and $m-1+
|J_1|\le m(1+|J_1|)$, the definition of the  expectation
$\E_{-,\ga_N}$,  the property \textbf{P}$(N-1,\ga_0,B)$,
(\ref{ecardJ}) and Lemma \ref{l725} imply
\begin{eqnarray*}
& &\E \left[ \1_{\Omega^c} \sum_{s \le \vert J_1 \vert}  \frac{
\vert \langle f_s \rangle_- \vert}{\langle {\bf Av} f \xi
\rangle_-} \right] = \E \left[ \1_{\Omega^c} \sum_{s \le \vert J_1
\vert} \E_{-,\ga_N}\lc \frac{ \vert \langle f_s \rangle_-
\vert}{\langle {\bf Av} f \xi
\rangle_-}\rc \right]\\
& &\quad \le \E \left[ \1_{\Omega^c} \sum_{s \le \vert J_1 \vert}
\frac{
(m-1+|J_1|)B\hat Q}{N-1} \right]\\
& &\quad \le 4 \frac{m}{N-1} B Q U_\infty \ \E \left[
\1_{\Omega^c} |J_1| (1+|J_1|)\ e^{4U_\infty R_1}
\right] \\
& &\quad \le 8 \frac{m}{N} B Q U_\infty\  \E \left[ \1_{\Omega^c}
|J_1| (1+|J_1|)\ e^{4U_\infty R_1} \right].
\end{eqnarray*}
Recall that, according to (\ref{easterisc}), we have
$$
\vert J_1 \vert \le \sum_{v \le R_1} \vert I_v^1 \vert,
$$
and that the quantity $R_1$ is a
${\rm Bin}( M,\frac{\ga}{N})$ random variable. Thus
\begin{eqnarray}
\E\lc \1_{\Omega^c} |J_1|\ e^{\la R_1}\rc&=&
 \E\lcl \E\lc \1_{\Omega^c}|J_1|\ e^{\la R_1}\Big|R_1 \rc\rcl=
\E\lcl e^{\la R_1}\ \E\lc \1_{\Omega^c} |J_1| \Big|R_1 \rc\rcl \nonumber \\
&\le &\ga \E\lc  R_1 e^{\la R_1}
\rc, \nonumber\\[3mm]
\E\lc \1_{\Omega^c}|J_1|^2 \ e^{\la R_1}\rc&=& \E\lcl
\E\lc\1_{\Omega^c} |J_1|^2\  e^{\la R_1}\Big|R_1
\rc\rcl \nonumber\\
&=&
\E\lcl e^{\la R_1}\ \E\lc \1_{\Omega^c}|J_1|^2 \Big|R_1 \rc\rcl
\nonumber \\
&\le& (\ga + \ga^2) \E\lc  ( R_1 + R_1^2) \ e^{\la
R_1} \rc.\label{nova2}
\end{eqnarray}
The proof of this proposition is now easily concluded by applying the
previous bounds, together with Lemma \ref{lbin}, to the quantity
$$\E \left[
\1_{\Omega^c} |J_1| (1+|J_1|)\ e^{4U_\infty R_1} \right].
$$

\hfill$\square$

\bigskip

We can turn  now to  the main aim of this Section:

\vspace{0.3cm}

 \noindent {\bf Proof of Theorem \ref{tproperty}:}
The result is now an immediate consequence of  (\ref{e712}) and
Proposition \ref{p726}, applied to
$$
B=B_0=\frac{ \al \ga^2 + \al^2 \ga^4}{1-\varepsilon},
$$ where
$\varepsilon$ satisfies
$$4 U_\infty \ \al \ga_0^2\ e^{4U_\infty}\
e^{\al \ga_0 (e^{4U_\infty}-1)}\ \left(3+ 2 \ga_0 + \al (\ga_0^2 +
\ga_0^3) e^{4U_\infty} \right)< \varepsilon <1.
$$

\hfill$\square$

\bigskip

Before closing this Section, we will give an easy consequence of Theorem
\ref{tproperty}: we will see that, as $N$ grows to $\infty$, the Gibbs measure
$G_N$ taken on a finite number of spins looks like a product measure. To this purpose, let us denote  by $\langle  \cdot  \rangle_\bullet$ the average with
respect to the product measure $\nu$ on $\Sigma_{N-1}$ such that
$$\forall i \le N-1,\quad \int \si_i\ d\nu(\mathop{\pmb{\rho}})
= \langle \si_i \rangle_-.$$ Equivalently, for a function $\bar f$
on $\Sigma_{N-1}$, we can write
$$\langle \bar f \rangle_\bullet=\langle \bar f(\si^1_1,\dots, \si_{N-1}^{N-1})\rangle_-,
$$
where $\si_i^i$ is the $i$-th coordinate of the $i$-th replica
$\mathop{\pmb{\rho}}^i$.
Recall also that, for $v\le R_1$, $I^1_v$ has been defined as
$$
I^1_v=\{i\le N-1;\ \ga_{i,k_v}=1\}.
$$
We now introduce the enumeration $\{i_1^v,\dots,i^v_{|I_v^1|}\}$
of this set. Furthermore, given the randomness contained in the
$\ga_{N,k}$, the law of $|I^1_v|$ is a ${\rm
Bin}(N-1,\frac{\ga}{N})$.
\begin{proposition}\label{p727}
Assume (\ref{e712}) and $\ga\le \ga_0$, and consider
$$\Theta=
\exp\lp \sum_{v\le R_1} \eta_{k_v} u\lp  \sum_{p\le |I^1_v|}
g_{i_p^v,k_v}\si_{i_p^v}+g_{N,k_v}\si_N\rp \rp . $$
Then, when $\Omega$
does not occur, we have
\begin{displaymath}
\E_{-,\ga_N}\left\vert \frac{\langle {\bf Av} \si_N \Theta
\rangle_-}{\langle {\bf Av}  \Theta \rangle_- } - \frac{\langle
{\bf Av} \si_N \Theta \rangle_\bullet}{\langle {\bf Av} \Theta
\rangle_\bullet } \right\vert \le 2 B_0 \ (|J_1|-1) \
\frac{|J_1|+1}{N-1} (e^{2 U_\infty}-1),
\end{displaymath}
where the conditional expectation $\E_{-,\ga_N}$ has been defined
at (\ref{emgn}).
\end{proposition}
\begin{remark}
The quantity $\Theta$ appears naturally in the decomposition of
the Hamiltonian
$H_{N,M}$. Indeed, on $\Omega_2^c$, we have
\begin{eqnarray*}
 &  &-H_{N,M} (\mathop{\pmb{\si}})\ =\ \sum_{k \le M} \eta_k u \left( \sum_{i \le N}
g_{i,k} \ga_{i,k} \si_{i} \right) \\
& = & \sum_{k \notin D^M_{N,1}} \eta_k u \left( \sum_{i \le N-1}
g_{i,k} \ga_{i,k} \si_{i} \right) + \sum_{k \in D^M_{N,1}} \eta_k
u \left( \sum_{i \le N-1} g_{i,k} \ga_{i,k} \si_{i} + g_{N,k}
\si_N
\right)\\
& = & \sum_{k \notin D^M_{N,1}} \eta_k u \left( \sum_{i \le N-1}
g_{i,k} \ga_{i,k} \si_{i} \right) +\sum_{v\le R_1} \eta_{k_v} u\lp
\sum_{p\le |I_v^1|} g_{i_p^v,k_v}\si_{i_p^v}+g_{N,k_v}\si_N\rp .
 \end{eqnarray*}
Observe also that $\xi$ defined by (\ref{defxi}) evaluated
for $n=1$ gives
$\xi=\Theta$.
\end{remark}

\vspace{0.3cm}

 \noindent {\bf Proof of
Proposition \ref{p727}}: The proof is similar to Proposition 7.2.7
in \cite{Tb}, and we include it here for sake of completeness:
On $\Omega^c$, since the sets $I^1_v$ are disjoint, the values
$i_p^v$, for any $v$ and $p$, are different and we can write
$$
\bigcup_{v\le R_1}I_v^1 =J_1\equiv \lcl j_1,\ldots j_{|J_1|} \rcl.
$$
Set
$$
f'=f'\big( \si_{j_1}^1,\ldots,\si_{j_{|J_1|}}^{1} \big) \equiv{\bf
Av} \si_N \Theta, \quad\quad
 f=f\big( \si_{j_1}^1,\ldots,\si_{j_{|J_1|}}^{1}
\big)\equiv{\bf Av}  \Theta.
$$
Let us also define, for $2\le l\le |J_1|$,
$$
f'_{j_l}=f'\big(\si_1^1,\ldots,\si_{j_l}^{j_l},
\si_{j_{l+1}}^{1},\ldots \si_{j_{|J_1|}}^{1}  \big)
$$ and $f_{j_l}$ in a similar way.
Then
\begin{align}
& \E_{-,\ga_N} \left\vert \frac{\langle {\bf Av} \si_N
\Theta \rangle_-}{\langle {\bf Av}  \Theta \rangle_- } -
\frac{\langle {\bf Av} \si_N \Theta \rangle_\bullet}{\langle {\bf
Av}  \Theta \rangle_\bullet} \right\vert\nonumber\\
& \quad = \E_{-,\ga_N} \left\vert
\frac{\langle f'_{j_1} \rangle_-}{\langle f_{j_1} \rangle_- } -
\frac{\langle f'_{j_{|J_1|}} \rangle_-}{\langle f_{j_{|J_1|}}
\rangle_- } \right\vert  \le  \sum_{2\le l\le |J_1|}
\E_{-,\ga_N}\lc \frac{\lla f'_{j_l-1}\rra_{-}}{\lla
f_{j_l-1}\rra_{-}} -\frac{\lla
f'_{j_l}\rra_{-}}{\lla f_{j_l}\rra_{-}}\rc\nonumber\\
&\quad\le \sum_{2\le l\le |J_1|}\left[ \E_{-,\ga_N}
\left\vert \frac{\langle f'_{j_l-1} - f'_{j_l} \rangle_- }{
\langle f_{j_l-1} \rangle_-} \right\vert + \E_{-,\ga_N} \left\vert
\frac{\langle f'_{j_l} \rangle_- \langle f_{j_l-1} - f_{j_l}
\rangle_- }{ \langle f_{j_l-1} \rangle_- - \langle f_{j_l}
\rangle_-} \right\vert\right].\label{e7.38}
\end{align}
Let us concentrate now on the first term of the right-hand side of (\ref{e7.38}),
since the other term can be bounded similarly: observe that, for
$2\le l\le |J_1|$, we have
$$
f'_{j_l}=f'_{j_{l-1}}\Delta, \quad\mbox{ with}\quad
e^{-2U_\infty}\le\Delta\le e^{2U_\infty}.
$$
Furthermore, it is easily seen that $f'_{j_l-1} - f'_{j_l}$ enjoys the
antisymmetric property assumed in Definition \ref{dproperty}.
Thus, applying $\textbf{P}(N-1,\ga_0,B_0)$, we get
$$
\E_{-,\ga_N} \left\vert \frac{\langle f'_{j_l-1} - f'_{j_l}
\rangle_- }{ \langle f_{j_l-1} \rangle_-} \right\vert \le
\frac{B_0(|J_1|+1)}{N-1} (e^{2 U_\infty}-1),$$ which ends the
proof.

\hfill$\square$


\section{Study of the magnetization}\label{magn}

For the non-diluted perceptron model, in the high temperature regime,
the asymptotic behavior of the magnetization can be summarized easily:
indeed, it has been shown in \cite{mcrt} that
$\langle \si_1\rangle$ converges in $L^2$ to a random variable of
the form $\tanh^2(z \sqrt{r})$, where $r$
is a solution to a deterministic equation, and $z\sim N(0,1)$.
Our goal in this section is to analyze
the same problem for the diluted perceptron model. However, in the current situation, the
limiting law is a more complicated object, and in order to present our asymptotic result, we
will go through a series of notations and preliminary lemmas.

\vspace{0.3cm}

Let ${\bf P}$ be the set of probability measures on $[-1,1]$. We
start by constructing a map $T:{\bf P}\to{\bf P}$ in the following
way: for any integer $\theta\ge 1$, let
$(\tau_1,\dots,\tau_\theta)$ be $\theta$ arbitrary integers. Then,
for $k=1,\dots,\theta$, let $t_k$ be the cumulative sum of the
$\tau_k$; that is, $t_0=0$ and $t_k=\sum_{\hat k\le k}
\tau_{\hat k}$ for $k\ge 1$. Let also $\{\bar g_{i,k},\ i,k \ge 1\}$
and $\{\bar g_{k},\ k \ge 1\}$ be two independent families of
independent standard Gaussian random variables. Define then a
random variable $\xi_{\theta,\tau}$ by
\begin{eqnarray}
\xi_{\theta,\tau}&=&
\xi_{\theta,\tau}(\si_1,\dots,\si_{t_{\theta}},\varepsilon)
\nonumber\\[2mm] &=&\exp \sum_{k=1}^\theta\
 u\left(\sum_{i=1}^{\tau_{k}} \bar g_{i,k}\ \si_{t_{k-1}+i}\ +
\bar g_k \ \varepsilon\right).\label{e740}
\end{eqnarray}
Whenever $\theta=0$, set also $\xi_\theta=1$, which is
equivalent to the convention $\sum_{k=1}^{0}w_k=0$ for any real sequence
$\{ w_k;k\ge 0 \}$.

\vspace{0.3cm}

Consider now $\mathop{\pmb{x}}=(x_1,\dots,x_{\sum_{k=1}^\theta \tau_k})$
with $|x_i|\le 1$ and a function
$$
f:\{-1,1\}^{\sum_{k=1}^\theta
\tau_k}\rightarrow \mathbb{R}.
$$
We denote  by  $\langle
f\rangle_{\mathop{\pmb{x}}}$ the average of $f$ with respect to
the product measure $\nu$ on $\{-1,1\}^{\sum_{k=1}^\theta \tau_k}$
such that $\int \si_i d\nu(\mathop{\pmb{\delta}})=x_i$, where
$\mathop{\pmb{\delta}}=(\si_1,\dots,\si_{\sum_{k=1}^\theta
\tau_k})$.
Using this notation, when $\theta \ge 1$, we define
$T_{\theta,\tau}:{\bf P}\to{\bf P}$ such that, for
$\mu\in{\bf P}$,
$T_{\theta,\tau}(\mu)$ is the law of
the random variable
\begin{equation}\label{e741}
\frac{\langle {\bf Av} \varepsilon \xi_{\theta,\tau}
\rangle_{\mathop{\pmb{X}}}}{\langle {\bf Av} \xi_{\theta,\tau}
\rangle_{\mathop{\pmb{X}}}},
\end{equation}
where $\mathop{\pmb{X}}=(X_1,\dots,X_{\sum_{k=1}^\theta \tau_k})$
is a sequence of i.i.d. random variables of law $\mu$ independent
of the randomness in $\xi_{\theta,\tau}$ and ${\bf Av}$ denotes
the average over $\varepsilon=\pm 1$. When $\theta=0$ we define
$T_{\theta,\tau}(\mu)$ as the Dirac measure at point 0.

Eventually, we can define the map $T:{\bf P}\to{\bf P}$ by
\begin{equation}\label{e7.42}
T(\mu)= \sum_{\theta\ge 0}\sum_{\tau_1,\dots,\tau_\theta\ge 0} \
\kappa(\theta,\tau_1,\dots,\tau_\theta)\
T_{\theta,\tau}(\mu),
\end{equation}
with
\begin{equation}\label{defkappa}
\kappa(\theta,\tau_1,\dots,\tau_\theta)= e^{-\alpha \gamma}\
\frac{(\alpha \gamma)^\theta}{\theta!}\ e^{-\theta \gamma}\
\frac{\gamma^{\sum_{l\le \theta} \tau_l}}{\tau_1!\cdots
\tau_\theta!},
\end{equation}
and where the coefficients $\al,\ga$ are the parameters of our
perceptron model. We will see that the asymptotic law $\mu$ of the magnetization
$\langle \si_1\rangle$ will satisfy the relation $\mu=T(\mu)$.
Hence, a first natural aim of this section is to prove that
the equation $\mu=T(\mu)$ admits a unique solution:
\begin{theorem}\label{t731}
Assume \begin{equation}\label{e750}2 U_\infty\ e^{2U_\infty}\
\alpha \gamma^2<\frac12.\end{equation}
 Then there exists a unique probability distribution
$\mu$ on $[-1,1]$ such that $\mu=T (\mu)$.
\end{theorem}
\begin{remark} Notice that (\ref{e712}) implies  (\ref{e750}).\end{remark}

In order to settle the fixed point argument for the proof of Theorem
\ref{t731}, we will need a metric on ${\bf P}$, and in fact it will be suitable for computational purposes to choose the
Monge-Kantorovich transportation-cost distance for the compact
metric space $([-1,1], | \cdot |)$: for two probabilities
$\mu_1$ and $\mu_2$ on $[-1,1]$, the distance between $\mu_1$ and $\mu_2$
will be defined as
\begin{equation}\label{defmong}
d(\mu_1,\mu_2)=\inf \E |X_1-X_2|,
\end{equation}
where this infimum is taken over all the pairs $(X_1,X_2)$ of
random variables such that the law of $X_j$ is $\mu_j$, $j=1,2$.
This definition is equivalent to say that $$d(\mu_1,\mu_2)=\inf
\int d(x_1,x_2) d\zeta(x_1,x_2), \quad\mbox{ with }\quad
d(x_1,x_2)=|x_2-x_1|,
$$
where this infimum is now taken over all probabilities $\zeta$ on
$[-1,1]^2$ with marginals  $\mu_1$ and $\mu_2$ (see Section 7.3 in
\cite{Tb} for more information about transportation-cost
distances). Finally, throughout this section, we also use
a local definition
of distance between two probabilities, with
respect to an event $\Omega$:
\begin{equation}\label{locdis}
d_\Omega(\mu_1,\mu_2)=\inf \E \left|(X_1-X_2)\1_\Omega\right|,
\end{equation}
where this infimum is as in (\ref{defmong}).

 \vspace{0.5cm}

\noindent {\bf Proof of Theorem \ref{t731}}: Assume that
$\theta\ge 1$ and $\tau_k\ge 1$ for some $k=1,\cdots,\theta$.
Then, using similar arguments to Lemma 7.3.5 in \cite{Tb} we can
prove, for $1\le i\le \sum_{k=1}^\theta \tau_k$, that
\begin{equation}\label{l735}
\left|\frac{\partial}{\partial x_i}\ \frac{\langle {\bf Av}
\varepsilon \xi_{\theta,\tau} \rangle_{\mathop{\pmb{x}}}}{\langle
{\bf Av} \xi_{\theta,\tau} \rangle_{\mathop{\pmb{x}}}}\right| \le
2U_\infty\ e^{2U_\infty},\end{equation} with
$\mathop{\pmb{x}}=(x_1,\dots,x_{\sum_{k=1}^\theta \tau_k})$. Then
if $\mathop{\pmb{y}}=(y_1,\dots,y_{\sum_{k=1}^\theta \tau_k})$,
the bound (\ref{l735}) implies that
\begin{equation}\label{e749}
\left|\frac{\langle {\bf Av} \varepsilon \xi_{\theta,\tau}
\rangle_{\mathop{\pmb{x}}}}{\langle {\bf Av} \xi_{\theta,\tau}
\rangle_{\mathop{\pmb{x}}}} - \frac{\langle {\bf Av} \varepsilon
\xi_{\theta,\tau} \rangle_{\mathop{\pmb{y}}}}{\langle {\bf Av}
\xi_{\theta,\tau} \rangle_{\mathop{\pmb{y}}}}\right| \le
2U_\infty\ e^{2U_\infty}\ \sum_{k=1}^\theta \sum_{i=1}^{\tau_k}
|x_{t_k+i}- y_{t_k+i}|.\end{equation} Remark that if $\theta=0$ or
$\theta \neq 0$ but $\tau_k=0$ for any $k=1,\dots,\theta$, then
the left-hand side of (\ref{e749}) is zero.

\vspace{0.3cm}

Let now $(X,Y)$ be a pair of random variables such that the laws of
$X$ and $Y$ are $\mu_1$ and $\mu_2$, respectively ($\mu_1$ and
$\mu_2$ are independent of the randomness in $\xi_{\theta,\tau}$).
Consider independent copies $(X_i,Y_i)_{i\le \sum_{k=1}^\theta
\tau_k}$ of this couple of random variables.
Then, if $\mathop{\pmb{X}}=(X_i)_{i\le
\sum_{k=1}^\theta \tau_k}$ and $\mathop{\pmb{Y}}=(Y_i)_{i\le
\sum_{k=1}^\theta \tau_k}$, we have that
\begin{displaymath} \frac{\langle {\bf Av} \varepsilon
\xi_{\theta,\tau} \rangle_{\mathop{\pmb{X}}}}{\langle {\bf Av}
\xi_{\theta,\tau} \rangle_{\mathop{\pmb{X}}}}
\stackrel{(d)}{=}
T_{\theta,\tau}(\mu_1)\quad {\rm and} \quad \frac{\langle {\bf Av}
\varepsilon \xi_{\theta,\tau} \rangle_{\mathop{\pmb{Y}}}}{\langle
{\bf Av} \xi_{\theta,\tau} \rangle_{\mathop{\pmb{Y}}}}
\stackrel{(d)}{=}
T_{\theta,\tau}(\mu_2).
\end{displaymath}
Hence, applying (\ref{e749}) for
$\mathop{\pmb{x}}=\mathop{\pmb{X}}$ and
$\mathop{\pmb{y}}=\mathop{\pmb{Y}}$ and taking first expectation
and then infimum over the choice of $(X,Y)$, we obtain
\begin{equation}\label{e748}
d(T_{\theta,\tau}(\mu_1),T_{\theta,\tau}(\mu_2))\le 2 U_\infty\
e^{2U_\infty}\ d(\mu_1,\mu_2) \sum_{k=1}^\theta \tau_k
.\end{equation}
Eventually, recall (see \cite[Lemma 7.3.2]{Tb}) that for a given sequence
 $\{ c_n;n\ge 1  \}$ of positive numbers such that $\sum_{n\ge 1} c_n=1$,
and two sequences $\{\mu_n,\nu_n;n\ge 1  \}$ of elements of ${\bf P}$, we have
\begin{equation}\label{sumsr}
d\lp \sum_{n\ge 1}c_n\mu_n,\,\sum_{n\ge 1}c_n\nu_n \rp
\le
\sum_{n\ge 1}c_n\, d\lp \mu_n,\nu_n\rp.
\end{equation}
Applying this elementary result to
$c_{\theta,\tau}=\kappa(\theta,\tau_1,\dots,\tau_\theta)$,
$\mu_{\theta,\tau}=T_{\theta,\tau}(\mu_1)$ and
$\nu_{\theta,\tau}=T_{\theta,\tau}(\mu_2)$, we get
\begin{align*}
&d(T(\mu_1),T(\mu_2))\le \sum_{\theta\ge
0}\sum_{\tau_1,\dots,\tau_\theta\ge 0} \
\kappa(\theta,\tau_1,\dots,\tau_\theta)\
d(T_{\theta,\tau}(\mu_1),T_{\theta,\tau}(\mu_2))
\\ & \le  2 U_\infty\ e^{2U_\infty}\
\sum_{\theta\ge 0}
\sum_{\tau_1,\dots,\tau_\theta\ge 0}
\sum_{k=1}^\theta
\kappa(\theta,\tau_1,\dots,\tau_\theta)\,
\tau_k\, d(\mu_1,\mu_2)\\
&= 2 U_\infty\ e^{2U_\infty}\ \left(\sum_{\theta\ge 0}\
e^{-\alpha \gamma}\ \frac{(\alpha \gamma)^\theta}{\theta !}\
\theta
\gamma\right)\ d(\mu_1,\mu_2)\\
&= 2 U_\infty\ e^{2U_\infty}\ \alpha \gamma^2 \ d(\mu_1,\mu_2),
\end{align*}
where we have used the fact  that the mean of a Poisson random variable
with parameter $\rho$ is $\rho$.
Then, under assumption (\ref{e750}), $T$ is a contraction and there
 exists a unique probability
distribution such that $\mu=T(\mu)$.

\hfill$\square$

\bigskip

Notice that the solution to the equation $\mu=T(\mu)$ depends on
the parameters $\alpha$ and $\gamma$. Furthermore, in the sequel,
we will need some continuity properties for the application
$(\al,\gamma)\mapsto\mu_{\al,\ga}$. Thus, we will set
$\mu=\mu_{\al,\ga}$ when we want to stress the dependence on the
parameters $\al$ and $\gamma$, and the following holds true:
\begin{lemma}\label{l737}
If $(\al,\gamma)$ and $(\al',\gamma')$ satisfy (\ref{e750}), then
$$d(\mu_{\al,\ga}, \mu_{\al',\gamma'})\le 4\left[|\gamma-\gamma'|
\alpha' \gamma' e^{|\gamma-\gamma'|}+ |\al \gamma-\al'\gamma'|
e^{|\al \gamma-\al'\gamma'|}\right].$$
\end{lemma}

\vspace{0.5cm}

\noindent {\bf Proof:} Since
$\mu_{\al,\ga}=T_{\al,\ga}(\mu_{\al,\ga})$ and
$\mu_{\al',\gamma'}=T_{\al',\gamma'}(\mu_{\al',\gamma'})$, using
the triangular inequality  and Theorem \ref{t731} we have
\begin{eqnarray*}d(\mu_{\al,\ga}, \mu_{\al',\gamma'})&\le &
d(T_{\al,\ga}(\mu_{\al,\ga}),T_{\al,\ga}(\mu_{\al',\gamma'}))
 + d(T_{\al,\ga}(\mu_{\al',\gamma'}),T_{\al',\gamma'}(\mu_{\al',\gamma'}))\\
&\le & \frac12 d(\mu_{\al,\ga},\mu_{\al',\gamma'}) +
d(T_{\al,\ga}(\mu_{\al',\gamma'}),T_{\al',\gamma'}(\mu_{\al',\gamma'})).\end{eqnarray*}
So
$$d(\mu_{\al,\ga}, \mu_{\al',\gamma'})\le 2 d(T_{\al,\ga}(\mu_{\al',\gamma'}),
T_{\al',\gamma'}(\mu_{\al',\gamma'}))$$ and we only need to deal
with
$d(T_{\al,\ga}(\mu_{\al',\gamma'}),T_{\al',\gamma'}(\mu_{\al',\gamma'}))$.
However, Lemma 7.3.3 in \cite{Tb} implies that
\begin{displaymath}\begin{array}{l}
\displaystyle
d(T_{\al,\ga}(\mu_{\al',\gamma'}),T_{\al',\gamma'}(\mu_{\al',\gamma'}))\\[2mm]
\qquad\qquad \displaystyle \le  2 \sum_{\theta\ge
0}\sum_{\tau_1,\dots,\tau_\theta\ge
0}\left|\kappa_{\al,\ga}(\theta,\tau_1,\dots,\tau_\theta)-
\kappa_{\al',\gamma'}(\theta,\tau_1,\dots,\tau_\theta)\right|\\[6mm]
\qquad\qquad \displaystyle \le  2 ( V_1 +
V_2),\end{array}\end{displaymath} with $\kappa$ defined in
(\ref{defkappa}) and
\begin{eqnarray*} V_1&=&\sum_{\theta\ge
0}\sum_{\tau_1,\dots,\tau_\theta\ge 0} e^{-\theta \gamma}
\frac{\gamma^{\sum_{l\le \theta} \tau_l}}{\theta!\ \tau_1!\cdots
\tau_\theta!}\left| e^{-\alpha \gamma}
(\alpha \gamma)^\theta- e^{-\alpha' \gamma'} (\alpha' \gamma')^\theta\right|,\\
 V_2&=&\sum_{\theta\ge
0}\sum_{\tau_1,\dots,\tau_\theta\ge 0} e^{-\alpha' \gamma'}
\frac{(\alpha' \gamma')^\theta}{\theta!\ \tau_1!\cdots
\tau_\theta!}
 \left|e^{-\theta \gamma}\gamma^{\sum_{l\le \theta} \tau_l}- e^{-\theta
 \gamma'}\gamma'^{\sum_{l\le \theta} \tau_l}
\right|.\end{eqnarray*}
Now, following the arguments of (7.53) in \cite{Tb},
we get
\begin{eqnarray*}
V_1&=&\sum_{\theta\ge 0} \frac{1}{\theta!} \left| e^{-\alpha
\gamma} (\alpha \gamma)^\theta- e^{-\alpha' \gamma'} (\alpha'
\gamma')^\theta\right|\le |\al \gamma-\al' \gamma'|
e^{|\al \gamma-\al'\gamma'|},\\
V_2&\le & |\gamma-\gamma'| e^{|\gamma-\gamma'|} \sum_{\theta\ge 0}
\theta e^{-\alpha' \gamma'} \frac{(\alpha'
\gamma')^\theta}{\theta!}=\alpha' \gamma' |\gamma-\gamma'|
e^{|\gamma-\gamma'|},\end{eqnarray*} which ends the proof of this
lemma.

\hfill$\square$

\bigskip

From now on, we will specialize our Hamiltonian to the case of interest
for us:
\begin{hypothesis}
The parameters  $\eta_k$, $k=1,\dots,M$ in the Hamiltonian
(\ref{ehamin}) are all equal to one.
\end{hypothesis}
This assumption  being made, we can now turn to the main result of the section:
\begin{theorem}\label{t741}
Let $\gamma_0$ be a positive number such that
\begin{equation}\label{oblidada}
 4 U_\infty \ \al \ga_0^2\ e^{4U_\infty}\ e^{\al \ga_0
(e^{4U_\infty}-1)}\ \left(3+ 2 \ga_0 + \al (\ga_0^2 + \ga_0^3)
e^{4U_\infty} \right) <1,
\end{equation}
and assume that there exists a positive number $C_0$ satisfying
\begin{equation}\label{e754}
C_0 \al \ \ga_0^6\ U_\infty\ e^{2U_\infty}   \le 1.
\end{equation}
Then for any $\gamma \le \gamma_0$, given any integer $m$, we can
find i.i.d. random variables $z_1,\dots,$ $z_m$ with law
$\mu_{\al,\ga}$ such that
\begin{equation}\label{limmag}
\E \lc \sum_{i\le m} \left\vert \langle
\si_i\rangle-z_i\right\vert \rc \le \frac{K m^3}{N},
\end{equation}
for a constant  $K>0$ independent of $m$.
\end{theorem}
\begin{remark}
The two conditions in the above theorem are met when the
following hypothesis is satisfied:  there exists $L>0$ such that
\begin{equation*}
L\ U_\infty\ \alpha \ \gamma_0^6\ \exp\left\{8U_\infty +\al \ga_0
\left(e^{4U_\infty}-1\right)\right\} \ <\ 1.
\end{equation*}
\end{remark}

As in the case of Theorem \ref{tproperty}, the proof of
Theorem \ref{t741} will require the introduction of some notations and preliminary Lemmas. Let us first recast relation (\ref{limmag}) in a suitable way for an induction procedure:
consider the metric space $[-1,1]^m$, equipped with the distance given by
$$ d((x_i)_{i\le m}, (y_i)_{i\le m})=  \sum_{i\le m} |x_i-y_i|.$$
We also denote by $d$ the transportation-cost distance on the
space of probability measures on $[-1,1]^m$, defined as in
(\ref{defmong}). Define now
\begin{equation}\label{e757}
D(N,M,m,\ga_0)=\sup_{\ga\le \ga_0} d\left(\mathcal{L}(\langle
\si_1\rangle,\dots,\langle \si_m\rangle),\mu_{\al,\ga}^{\otimes
m}\right),
\end{equation}
where $\mathcal{L}(X)$ stands for the law of the
random variable $X$. Then the statement of Theorem \ref{t741} is equivalent to say that, under Hypothesis (\ref{e754}), we have
$$D(N,M,m,\ga_0)\le \frac{Km^3}{N},$$ for any fixed integer $m\ge 1$.

\vspace{0.3cm}

It will also be useful to introduce a cavity formula for $m$
spins, which we proceed to do now: generalizing  some aspects of
the previous section, we consider, for $p\in \{1,\dots,m\}$, the
random sets
$$D_{N,p}^M=\{k\le M;\ \ga_{N-p+1,k}=1\},$$
and
$$F_{N,M}^m=\bigcup_{p=1}^m D_{N,p}^M\ .$$
We also define the following two rare events:
\begin{eqnarray*}
\tilde \Omega_1&=&\{\exists k\le M,\ p_1, p_2 \le m;\,
\ga_{N-p_1+1,k}=\ga_{N-p_2+1,k}=1\},\\
\tilde \Omega_2&=&\{\exists i\le N-m,\ k_1, k_2 \in F_{N,M}^m;\,
\ga_{i,k_1}=\ga_{i,k_2}=1\},
\end{eqnarray*}
satisfying
\begin{equation}\label{bndptan}
P(\tilde \Omega_1)\le \frac{\al \ga^2 m^2}{N}\qquad {\rm and}
\qquad P(\tilde \Omega_2)\le \frac{\al^2 \ga^4
m^2}{N}.
\end{equation}
Then, the following properties hold true: first, for a fixed $k$, if
$\tilde\Omega_1^c$ is realized, we have
$$
{\rm Card} \{p\le m,\ \ga_{N-p+1,k}=1\}\le 1.
$$
Moreover, still on  $\tilde \Omega_1^c$, for $p_1\neq p_2$,
$$D_{N,p_1}^M \cap D_{N,p_2}^M=\emptyset;$$
and hence,
$$
R_m\equiv |F_{N,M}^m|= \sum_{p=1}^m |D_{N,p}^M| =\sum_{k\le M} \sum_{p\le
m} \ga_{N-p+1,k}.$$ Actually, notice that we always have
$$R_m \le \sum_{p=1}^m |D_{N,p}^M|.$$
Let us introduce now an enumeration of $F_{N,M}^m$:
$$F_{N,M}^m=\{k_1,\dots,k_{R_m}\},$$
and for any $v \le R_m$ set
$$I_v^m=\{j\le N-m; \ga_{j,k_v}=1\}.$$
Then, on $\tilde \Omega_2^c$, we get
\begin{equation}\label{einterbuiad}I_{v_1}^m\cap
I_{v_2}^m=\emptyset,\quad {\rm if} \ v_1\neq v_2,
\end{equation}
and we can also write
$$J_m=\bigcup_{v\le R_m} I_v^m
=\bigcup_{v\le R_m} \{ j\le N-m;\ \ga_{j,k_v}=1 \}.$$

\vspace{0.3cm}

Let us separate now the $m$ last spins in the Hamiltonian
$H_{N,M}$: if $\tilde\Omega_1^c$ is realized, for
$\mathop{\pmb{\rho}}= (\si_1,\dots,\si_{N-m})$,
we have the following decomposition:
\begin{equation*}
-H_{N,M}(\mathop{\pmb{\si}}) = -H_{N-m,M}^-(\mathop{\pmb{\rho}}) +
\log \xi,
\end{equation*}
with
\begin{align}
&-H_{N-m,M}^-(\mathop{\pmb{\rho}})= \sum_{k \in (F_{N,M}^m)^c}
  u \left( \sum_{i \le N-m} g_{i,k} \ga_{i,k} \si_{i}
\right),\nonumber\\[3mm]
&\xi = \exp\lp\sum_{p=1}^m \sum_{k\in D_{N,p}^M}
u \left( \sum_{i
\le N-m} g_{i,k} \ga_{i,k} \si_{i} + g_{N-p+1,k}
\si_{N-p+1}\right)\rp.\label{defxim}
\end{align}
Observe that, in the last formula, $H_{N-m,M}^-(\mathop{\pmb{\rho}})$ is not exactly the
Hamiltonian of a $(N-m)$-spin system  changing $\ga$ into
$\ga^{-}$, because the set $F_{N,M}^n$ is not deterministic.
But this problem
will be solved again by conditioning  upon the random variables
$\{ \ga_{N-p+1,k},\
p=1,\dots,m, \ k\le M \}$. For the moment, let us just mention that the $m$ cavity formula will be the following:
given $f$ on $\Sigma_N$, we have
\begin{equation}\label{nova1}
\1_{\tilde \Omega_1^c}\langle f \rangle = \1_{\tilde
\Omega_1^c}\frac{\langle {\bf Av} f \xi \rangle_{-}}{\langle {\bf
Av} \xi \rangle_{-}},
\end{equation}
 where $\langle \cdot
\rangle_-$ is the average with respect to $H_{N-m,M}^-$ and ${\bf
Av}$ is the average with respect to last $m$ spins. Moreover, in
the last formula, we have kept the notation $\xi$ from Section
\ref{spincor}, which hopefully will not lead to any confusion.
Eventually, we denote by  $\mathcal{L}_0$ the law of a random
variable conditioned by $\{ \ga_{N-p+1,k},\ p=1,\dots,m, \ k\le M
\}$, and by $\E_{\ga_{N,m}}$ the associated conditional
 expectation.

\vspace{0.3cm}

We can start now stating and proving the lemmas and propositions that
will lead to the proof of Theorem \ref{t741}. Recall that
given $\mathop{\pmb{x}}=(x_1,\dots,x_{N-m})$, $|x_i|\le 1$, and a
function $f$ on $\Sigma_{N-m}$, $\langle f
\rangle_{\mathop{\pmb{x}}}$ means the average of $f$ with respect
to the product measure $\nu$ on $\Sigma_{N-m}$ such that $\int
\si_i d\nu(\mathop{\pmb{\rho}})=x_i$, for $0\le i \le N-m$. Recall also
that $\ga^-=\ga \frac{N-m}{N}$.
Then, as a direct consequence of the definition of the operator
$T_{\theta,\tau}$, we have the following result:
\begin{lemma}\label{l744}
Let $\mathop{\pmb{X}}=(X_1,\dots,X_{N-m})$ be an independent
sequence of random variables, where the law of each $X_l$ is
$\mu_{\al,\ga^-}$. Set
$$w_p=\frac{\langle {\bf Av} \si_{N-p+1} \xi
\rangle_{\mathop{\pmb{X}}}}{\langle {\bf Av} \xi
\rangle_{\mathop{\pmb{X}}}},\qquad p=1,\dots,m.$$ Then, on
$\tilde\Omega^c=(\tilde \Omega_1\cup \tilde \Omega_2)^c$, we have
\begin{displaymath}
\mathcal{L}_0(w_1\dots,w_m)=T _{|D_{N,1}^M|,(|I_k^m|, k\in
D_{N,1}^M)}(\mu_{\al,\ga^-})\otimes \cdots \otimes
T_{|D_{N,m}^M|,(|I_k^m|, k\in D_{N,m}^M)}(\mu_{\al,\ga^-}).
\end{displaymath}
\end{lemma}

We will now try to relate the random variables $w_p$ with the
magnetization of the $m$ last spins. A first step in that
direction is the following lemma where we use the random value of
the parameter $\al^-$ associated to the Hamiltonian of a
$(N-m)$-spin system.

\begin{lemma}\label{l744bis}
On $\tilde\Omega^c$, set
\begin{displaymath}\Gamma_m=d\left(\mathcal{L}_0(w_1\dots,w_m),
\mathcal{L}_0(\bar w_1\dots,\bar w_m)\right),\end{displaymath}
where, for $p=1,\dots,m$,
$$\bar w_p=T _{|D_{N,p}^M|,(|I_k^m|, k\in
D_{N,p}^M)}(\mu_{\al^-,\ga^-}),\quad {\rm with} \
\al^-=\frac{M-R_m}{N-m}.$$ Then, on $\tilde\Omega^c$, we have
\begin{displaymath}
\Gamma_m\le 2 U_\infty\ e^{2U_\infty}\ \ga_0\ \frac{|R_m-m
\al|}{N} \exp\left\{\ga_0\ \frac{|R_m-m \al|}{N}\right\}\
\sum_{k=1}^{R_m} |I_k^m|.
\end{displaymath}
\end{lemma}

\vspace{0.5cm}

 \noindent {\bf Proof:}
Using (\ref{e748}) we obtain
 \begin{equation*}
\Gamma_m\le  2 U_\infty\ e^{2U_\infty}\ \sum_{p=1}^m \sum_{k\in
D_{N,p}^M}|I_k^m|\
d(\mu_{\al,\ga^-},\mu_{\al^-,\ga^-}).\end{equation*}
The proof of
this lemma is then easily finished thanks to Lemma
 \ref{l737}, and taking the following equality  into
account:
 \begin{equation*}\label{ealfa}\ga^-\ |\al-\al^-|= \ga
 \frac{|R_m-m\al|}{N}.\end{equation*}

  \hfill$\square$

\bigskip

Notice that we have introduced the random variables $\bar w_p$
for the following reason:
given the randomness contained in the $\{\ga_{N-p+1,k} ,\
p=1,\dots,m, k\le M\}$, $\bar w_p$ can be
interpreted as
$$\bar w_p=\frac{\langle {\bf Av} \si_{N-p+1} \xi
\rangle_{\mathop{\pmb{\bar X}}}}{\langle {\bf Av} \xi
\rangle_{\mathop{\pmb{\bar X}}}},\qquad p=1,\dots,m,$$ where
$\mathop{\pmb{\bar X}}=(\bar X_1,\dots,\bar X_{N-m})$ is an
independent sequence of random variables with law
$\mu_{\al^-,\ga^-}$.

\begin{lemma}\label{l743}
Consider $\mathop{\pmb{Z}}=(\langle \si_1\rangle_-,\dots,\langle
\si_{N-m}\rangle_-)$,  and denote
$$u_p=\frac{\langle {\bf Av} \si_{N-p+1} \xi
\rangle_{\mathop{\pmb{Z}}}}{\langle {\bf Av} \xi
\rangle_{\mathop{\pmb{Z}}}},\qquad p=1,\dots,m.$$ Then, on
$\tilde\Omega^c$,
\begin{multline*}
d\left(\mathcal{L}_0(u_1,\dots,u_m), \mathcal{L}_0(\bar w_1,\dots,\bar w_m)
\right)\\
\le 4 D(N-m, |(F_{N,M}^m)^c|, |J_m|,\ga_0)\ U_\infty \
e^{2U_\infty},
\end{multline*}
where the quantity $D$ has been defined at relation (\ref{e757}).
\end{lemma}

\vspace{0.5cm}

 \noindent {\bf Proof:}  As in (\ref{l735}) we can obtain, for any
$i \le N-m$,
\begin{equation}\label{e768}\left|\frac{\partial}{\partial x_i}\ \frac{\langle {\bf Av}
\si_{N-p+1} \xi \rangle_{\mathop{\pmb{x}}}}{\langle {\bf Av} \xi
\rangle_{\mathop{\pmb{x}}}}\right| \le 2U_\infty \
e^{2U_\infty}.\end{equation}
But in fact, these derivatives are vanishing, unless
$$i\in I_{N,p}^m \equiv \bigcup_{v;\ k_v\in D_{N,p}^M} I^m_v,$$
for some $p=1,\dots,m$. Indeed, on $\tilde\Omega^c$, from
(\ref{einterbuiad}), we have
$$I_{N,p_1}^m \cap I_{N,p_2}^m = \emptyset,\qquad {\rm if}\
p_1\neq p_2.$$ Then, for a given $p\in \{1,\dots,m\}$, we can
decompose $\xi$ into $\xi=\xi_{N,p}\ \bar \xi_{N,p}$, with
\begin{eqnarray*}
\xi_{N,p}&=& \exp\sum_{k\in D_{N,p}^M}   u \left( \sum_{i \le N-m}
g_{i,k} \ga_{i,k} \si_{i} + g_{N-p+1,k}
\si_{N-p+1}\right)\\[3mm]
&=& \xi_{N,p}\left(\{\si_i, i\in I_{N,p}^m\},
\si_{N-p+1}\right),\\[3mm] \bar \xi_{N,p}&=&\bar
\xi_{N,p}\left(\{\si_i, i\in J_m\backslash I_{N,p}^m\},
\si_{N-\bar p+1}, \bar p\le m, \bar p\neq p\}\right).
\end{eqnarray*}
Then \begin{displaymath} \frac{\langle {\bf Av} \si_{N-p+1} \xi
\rangle_{\mathop{\pmb{x}}}}{\langle {\bf Av} \xi
\rangle_{\mathop{\pmb{x}}}}=\frac{\langle {\bf Av} \si_{N-p+1}
\xi_{N,p} \rangle_{\mathop{\pmb{x}}}\ \langle {\bf Av} \bar
\xi_{N,p} \rangle_{\mathop{\pmb{x}}}}{\langle {\bf Av} \xi_{N,p}
\rangle_{\mathop{\pmb{x}}}\ \langle {\bf Av} \bar \xi_{N,p}
\rangle_{\mathop{\pmb{x}}}}=\frac{\langle {\bf Av} \si_{N-p+1}
\xi_{N,p} \rangle_{\mathop{\pmb{x}}}}{\langle {\bf Av} \xi_{N,p}
\rangle_{\mathop{\pmb{x}}}},
\end{displaymath}
and clearly the derivative $\frac{\partial}{\partial x_i}$ is zero
when $i$ does not belong to $I_{N,p}^m$, for any $p \in
\{1,\ldots,m\}$.

\vspace{0.3cm}

Now, invoking inequality (\ref{e768}), we get
$$ \sum_{p=1}^m \left| \frac{\langle {\bf Av} \si_{N-p+1} \xi
\rangle_{\mathop{\pmb{\bar X}}}}{\langle {\bf Av} \xi
\rangle_{\mathop{\pmb{\bar X}}}}-u_p\right|\le \left( \sum_{p=1}^m
\sum_{i\in I_{N,p}^m} \left|\bar X_i-\langle
\si_i\rangle_-\right|\right)\ 2 U_\infty \ e^{2U_\infty}.$$ Then,
the definition of $\E_{\ga_{N,m}}$ and (\ref{e757}) easily yield
$$\E_{\ga_{N,m}}\left( \sum_{p=1}^m \sum_{i\in
I_{N,p}^m} \left|\bar X_i-\langle \si_i\rangle_-\right|\right)\le
2 D(N-m, |(F_{N,M}^m)^c|, |J_m|,\ga_0),$$ which ends the proof.

  \hfill$\square$

\bigskip

Set now, for $1\le p \le m$,
\begin{equation}\label{deful}
\bar u_p=\frac{\langle {\bf Av} \si_{N-p+1} \xi \rangle_{-}}
{\langle {\bf Av} \xi \rangle_{-}} .
\end{equation}
Then $\bar u_p$ is closer to the real magnetization in the sense
that $\bar u_p=\langle \si_{N-p+1}\rangle$ on $\tilde \Omega^c$,
and the following Lemma claims that the distance between $\bar
u_p$ and $u_p$ vanishes as $N\to\infty$.
\begin{lemma}\label{l742} For $1 \le p \le m$, let $\bar u_p$ be defined by
(\ref{deful}). Then, on $\tilde\Omega^c$, we have
$$d\left(\mathcal{L}_0(\bar u_1,\dots,\bar u_m),
\mathcal{L}_0(u_1,\dots,u_m)\right)\le 2 B_0 \
\frac{|J_m|^2-1}{N-m+1}\ (e^{2 U_\infty}-1),$$ where the constant
$B_0$ has been defined in the previous section.\end{lemma}

\vspace{0.5cm}

\noindent {\bf Proof:}
The computations can be leaded here almost like in the proof of Proposition
\ref{p727}, and the details are left to the reader.

\hfill$\square$

\bigskip

We will now identify the law of the $\bar u_p$ in terms of laws of
the type $T(\mu_{\al,\ga^-})$:
\begin{lemma}\label{l745}
Recall that $d_{\tilde \Omega^c}$ has been defined by relation
(\ref{locdis}). Then, for $m\ge 1$, set
\begin{multline*}
\delta_m
=
 d_{\tilde \Omega^c}\Big(\mathcal{L}(\bar u_1,\dots,\bar u_m), \sum_{(b)}
\sum_{(\mathop{\pmb{v}})}
 a\left((b_1,\mathop{\pmb{v_1}}),\dots, (b_m,\mathop{\pmb{v_m}})\right)\\
T_{b_1,\mathop{\pmb{v_1}}}(\mu_{\al,\ga^-}) \otimes \cdots \otimes
T_{b_m,\mathop{\pmb{v_m}}}(\mu_{\al,\ga^-})\Big),
\end{multline*}
where we have used the following conventions: for $j\le m$,
$\mathop{\pmb{v_j}}$ is a multi-index of the form
$\mathop{\pmb{v_j}}=(v_1^j,\dots,v_{b_j}^j)$; the first summation
$\sum_{(b)}$ is over $b_j\ge 0$, for $j=1,\dots,m$; the second one
$\sum_{(\mathop{\pmb{v}})}$ is over $v_1^j,\dots,v_{b_j}^j\ge 0$,
for $j=1,\dots,m$; and $a((b_1,\mathop{\pmb{v_1}}),\dots,
(b_m,\mathop{\pmb{v_m}}))$ is defined by
$$
a\left((b_1,\mathop{\pmb{v_1}}),\dots,
(b_m,\mathop{\pmb{v_m}})\right) =P\left(|D_{N,j}^M|=b_j,(|I_k^m|,
k\in D_{N,j}^M)=\mathop{\pmb{v_j}}, \forall j\le m\right).
$$
Then, under the conditions of Lemma \ref{l742}, we have
$$
\delta_m\le c_1(N,m),
$$
with
\begin{eqnarray*}
c_1(N,m)&=& 4 \ U_\infty \ e^{2U_\infty}\
\E\left(D(N-m,|(F_{N,M}^m)^c|, |J_m|,\ga_0)\right)\\
&&+2 B_0 \ \frac{\E|J_m|^2-1}{N-m+1}\ (e^{2 U_\infty}-1)\\&& +2
U_\infty\ e^{2U_\infty}\  \frac{\ga_0}{N}  \E\left(|R_m-m \al| \
|J_m| \exp\left\{\frac{\ga_0}{N} |R_m-m \al|\right\} \right).
\end{eqnarray*}
\end{lemma}

\vspace{0.5cm}

\noindent {\bf Proof:} This result is easily obtained by combining
Lemmas \ref{l744}, \ref{l744bis}, \ref{l743}, \ref{l742} and
taking expectations.

 \hfill$\square$

\bigskip

With Lemma \ref{l745} in hand, we can  see that the remaining task left
 to us is mainly to compare the coefficients
$a((b_1,\mathop{\pmb{v_1}}),\dots, (b_m,\mathop{\pmb{v_m}}))$ with
the coefficients $\kappa_{\al,\ga^-}(b_j,\mathop{\pmb{v_j}})$.
This is done in the following lemma.
\begin{lemma}\label{l746}
With the conventions of Lemma \ref{l745},
we have
\begin{equation}\label{e775}\sum_{(b)}
\sum_{(\mathop{\pmb{v}})} \Bigg|
a\left((b_1,\mathop{\pmb{v_1}}),\dots,
(b_m,\mathop{\pmb{v_m}})\right)- \prod_{j=1}^m
\kappa_{\al,\ga^-}(b_j,\mathop{\pmb{v_j}})\Bigg| \le \frac{ m
L_0(\ga) }{N}.
\end{equation}
\end{lemma}

\vspace{0.5cm}

\noindent {\bf Proof:}
In fact, it is easily seen that we only need to prove that
\begin{displaymath}
\sum_{t,\mathop{\pmb{v}}\ge 0} \left|a(b,\mathop{\pmb{v}})-
\kappa_{\al,\ga^-}(b,\mathop{\pmb{v}})\right|\le \frac{L_0(\ga)
}{N},
\end{displaymath}
with $\mathop{\pmb{v}}=(v_1,\dots,v_b)$. However, notice that
\begin{displaymath}
a(b,\mathop{\pmb{v}})=\binom{M}{b}
\left(\frac{\gamma}{N}\right)^b
\left(1-\frac{\gamma}{N}\right)^{M-b} \prod_{l=1}^b
\binom{N-m}{v_l} \left(\frac{\gamma}{N}\right)^{v_l}
\left(1-\frac{\gamma}{N}\right)^{N-m-v_l},
\end{displaymath}
and recall that \begin{displaymath}
\kappa_{\al,\ga^-}(b,\mathop{\pmb{v}})=e^{-\alpha \gamma^{-}}
\frac{(\alpha \gamma^{-})^{b}}{b!}\ e^{-b \gamma^{-}}
\frac{(\gamma^{-})^{\sum_{l\le b}v_l}}{v_1!\cdots
v_b!}.\end{displaymath} Then
\begin{displaymath}
\sum_{b,\mathop{\pmb{v}}\ge 0}
\left|a(b,\mathop{\pmb{v}})-\kappa_{\al,\ga^-}(b,\mathop{\pmb{v}})
\right| \le A+B,
\end{displaymath}
with
\begin{eqnarray*}
A&=&\sum_{b,\mathop{\pmb{v}}\ge 0} \left|e^{-\alpha \gamma^{-}}
\frac{(\alpha \gamma^{-})^{b}}{b!} \ \bar A_{b,\mathop{\pmb{v}}}\right|,\\
\bar A_{b,\mathop{\pmb{v}}}&=& \left| e^{-b \gamma^{-}}
\frac{(\gamma^{-})^{\sum_{l\le b}v_l}}{v_1!\cdots
v_b!}-\prod_{l=1}^b \binom{N-m}{v_l}
\left(\frac{\gamma}{N}\right)^{v_l}
\left(1-\frac{\gamma}{N}\right)^{N-m-v_l} \right|,\\
B&=&\sum_{b,\mathop{\pmb{v}}\ge 0} \left| \bar B_b \prod_{l=1}^b
\binom{N-m}{v_l} \left(\frac{\gamma}{N}\right)^{v_l}
\left(1-\frac{\gamma}{N}\right)^{N-m-v_l}\right|,\\
\bar B_b&=& \left| e^{-\alpha \gamma^-} \frac{(\alpha
\gamma^-)^{b}}{b!}- \binom{M}{b} \left(\frac{\gamma}{N}\right)^b
\left(1-\frac{\gamma}{N}\right)^{M-b} \right|.
\end{eqnarray*}
Now, following the estimates for the approximation of a Poisson
distribution by a Binomial given in \cite[Lemma 7.4.6]{Tb}, we
can bound $\bar A_{b,\mathop{\pmb{v}}}$ and $\bar B_b$ by a
quantity of the form $\frac{c}{N}$. The proof is then easily
finished.

\hfill$\square$

\bigskip

Let us relate now the law of $(\bar u_1,\ldots,\bar u_m)$ with
$\mu_{\al,\ga^-}^{\otimes m}$.
\begin{lemma}\label{lea}
We have
$$
d_{\tilde \Omega ^c}\Big(\mathcal{L}(\bar u_1,\dots,\bar u_m),
\mu_{\al,\ga^-}^{\otimes m} \Big) \le c_2(N,m),
$$
with
\begin{eqnarray*}
c_2(N,m)&=& 4 \ U_\infty \ e^{2U_\infty}\
\E\left(D(N-m,|(F_{N,M}^m)^c|, |J_m|,\ga_0)\right)\\
&&+2 B_0 \ \frac{\E|J_m|^2-1}{N-m+1}\ (e^{2 U_\infty}-1)+\frac{2
m^2L_0(\ga_0)}{N}\\
&& +2 U_\infty\ e^{2U_\infty}\  \frac{\ga_0}{N}  \E\left(|R_m-m
\al| \ |J_m| \exp\left\{\frac{\ga_0}{N} |R_m-m \al|\right\}
\right).
\end{eqnarray*}
\end{lemma}

\vspace{0.5cm}

\noindent {\bf Proof:} Notice that, invoking relation
(\ref{e7.42}) and Theorem  \ref{t731}, we get
\begin{equation*}
\sum_{(b)} \sum_{(\mathop{\pmb{v}})} \Big( \prod_{j=1}^m
\kappa_{\al,\ga^-}(b_j,\mathop{\pmb{v_j}}) \Big)
T_{b_1,\mathop{\pmb{v_1}}}(\mu_{\al,\ga^-}) \otimes \cdots \otimes
T_{b_m,\mathop{\pmb{v_m}}}(\mu_{\al,\ga^-})=\mu_{\al,\ga^-}^{\otimes
m}.
 \end{equation*}
Then, the results follows easily from Lemmas \ref{l745} and
\ref{l746}, Lemma 7.3.3 in \cite{Tb}  and the triangular
inequality.

\hfill$\square$

\bigskip

We are now ready to end the proof of the main result concerning
the magnetization of the system.

\vspace{0.5cm}

\noindent {\bf Proof of Theorem \ref{t741}}: First of all, notice
that by symmetry we have
$$\mathcal{L}(\langle
\si_1\rangle,\dots,\langle \si_m\rangle)=\mathcal{L}(\langle
\si_{N-m+1}\rangle,\dots,\langle \si_N\rangle).$$ Furthermore,
thanks to (\ref{nova1}) and (\ref{bndptan}) and Lemma \ref{l737},
we can write
\begin{eqnarray*} & &D(N,M,m,\ga_0) =  \sup_{\ga\le \ga_0}
d\left(\mathcal{L}(\langle \si_1\rangle,\dots,\langle
\si_m\rangle),\mu_{\al,\ga}^{\otimes m}\right)\\
&  & \quad \le\sup_{\ga\le \ga_0} d_{\tilde \Omega
^c}\left(\mathcal{L}\left(\frac{\langle {\bf Av} \si_{N-m+1} \xi
\rangle_{-}}{\langle {\bf Av} \xi \rangle_{-}},\dots,
\frac{\langle {\bf Av} \si_{N} \xi \rangle_{-}}{\langle {\bf Av}
\xi \rangle_{-}} \right),\mu_{\al,\ga}^{\otimes m}\right)\\
&& \qquad +\frac{2m^3 \al \ga_0^2(1+\al \ga_0^2)}{N}
\\
&   & \quad \le \sup_{\ga\le \ga_0} d_{\tilde \Omega
^c}\left(\mathcal{L}\left(\frac{\langle {\bf Av} \si_{N-m+1} \xi
\rangle_{-}}{\langle {\bf Av} \xi \rangle_{-}},\dots,
\frac{\langle {\bf Av} \si_{N} \xi \rangle_{-}}{\langle {\bf Av}
\xi \rangle_{-}} \right),\mu_{\al,\ga^-}^{\otimes m}\right)\\
&& \qquad +\frac{2m^3 \al \ga_0^2(1+\al \ga_0^2)}{N}\\
&& \qquad + \frac{4m\al \ga_0}{N} \left( \ga_0
 \exp\Big\{\frac{m
\ga_0}{N}\Big\}+ \exp\Big\{\frac{m \al \ga_0}{N}\Big\} \right).
\end{eqnarray*}
 Then, Lemma \ref{lea} implies
 \begin{eqnarray*}
& &D(N,M,m,\ga_0) \le 4 \ U_\infty \ e^{2U_\infty}\
\E\left(D(N-m,|(F_{N,M}^m)^c|, |J_m|,\ga_0)\right) \\
&& \qquad +2 B_0 \ \frac{\E|J_m|^2-1}{N-m+1}\ (e^{2 U_\infty}-1)\\
&& \qquad +2 U_\infty\ e^{2U_\infty}\  \frac{\ga_0}{N}
\E\left(|R_m-m \al| \ |J_m| \exp\left\{\frac{\ga_0}{N} |R_m-m
\al|\right\}
\right)\\
&& \qquad +\frac{2 m^2L_0(\ga_0)}{N} + \frac{12 m^3 \al \ga_0^4
\exp(\ga_0)}{N}.
\end{eqnarray*}
It is readily checked, as we did in  (\ref{nova2}), that
\begin{eqnarray*}
\E \left( \left| J_m \right| \right)\  & \le & \frac{N-m}{N}\ \al
m\ga_0^2,\\
\E \left( \left| J_m \right|^2 \right)  & \le &  \frac{N-m}{N}
(\ga_0 +
\ga_0^2) \big( \al m \ga_0 + (\al m \ga_0 )^2 \big) , \\
\E \left( \left| J_m \right|^3 \right)  & \le &  \frac{N-m}{N}
(\ga_0 + 3\ga^2_0 + \ga^3_0) \big( \al m \ga_0 + 3(\al m \ga_0 )^2
+ (\al m \ga_0 )^3\big).
\end{eqnarray*}
Thus, using the fact that $R_m \le Y$ 
where $Y \sim B(mM, \frac{\gamma}{N})$, together with the trivial bound
$R_m \le M$, there exists a constant $K_0 \ge 1$ such that
\begin{eqnarray} D(N,M,m,\ga_0) &\le&
4 \ U_\infty \ e^{2U_\infty} \E\left(D(N-m,|(F_{N,M}^m)^c|,
|J_m|,\ga_0)\right)
\nonumber \\
& & +  \frac{K_0 m^3\big[\al \ga_0^4 \exp(\frac32\ga_0)+
L_0(\ga_0)\big] }{N}.\label{nova3}
\end{eqnarray}
Now we are able to prove, by induction over $N$, that
$$
D(N,M,m,\ga_0 ) \le  \frac{2 K_0 m^3\big[\al \ga_0^4
\exp(\frac32\ga_0)+ L_0(\ga_0)\big] }{N}, \qquad \mbox{for all}
\quad m \le \frac{N}{2}.$$ Indeed, in order to check the induction
step from $N-1$ to $N$, notice that $|(F_{N,M}^m)^c| \le M$ and
that
$$
\E \left( \left| J_m \right|^3 \right)   \le  25  \frac{N-m}{N}
(m^3 \al \ga_0^6).$$ So, using also that
$$
P \left( \left| J_m \right| \ge \frac{N}{2} \right) \le
\frac{4}{N^2} \E \left( \left| J_m \right|^2 \right) \le \frac{16
m^2 \al \ga^4_0}{N^2},$$
 and  by our induction
hypothesis and (\ref{nova3}), we have
\begin{eqnarray*} & &D(N,M,m,\ga_0) \le
 \frac{K_0 m^3\big[\al \ga_0^4 \exp(\frac32\ga_0)+
L_0(\ga_0)\big] }{N}
\nonumber \\
& & \ + 4 \ U_\infty \ e^{2U_\infty} \left( \frac{2K_0 \E
\left( \left| J_m \right|^3 \right)  \big[\frac{M}{N-m}\ \ga_0^4
\exp(\frac32\ga_0)+ L_0(\ga_0)\big]}{N-m} + \frac{32 m^3 \al
\ga_0^4}{N^2} \right)\\
& & \ \le \frac{K_0 m^3\big[\al \ga_0^4 \exp(\frac32\ga_0)+
L_0(\ga_0)\big] }{N} \\
& & \ + 4 \ U_\infty \ e^{2U_\infty} \left( \frac{50 K_0 m^3
\al \ga_0^6  \big[2\al \ \ga_0^4 \exp(\frac32\ga_0)+
L_0(\ga_0)\big]}{N} + \frac{32 m^3 \al \ga_0^4}{N^2} \right).
\end{eqnarray*}
Finally, since $M < N-m$, the proof easily follows from
hypothesis (\ref{e754}).

\hfill$\square$

\section{Replica symmetric formula}\label{symf}
Now that the limiting law of the magnetization has been computed,
we can try to evaluate the asymptotic behavior of the free energy
of our system, namely
$$
p_N(\gamma) = \frac{1}{N} \E \left[ \log \left(
\sum_{\mathop{\pmb{\si}}\in \Sigma_N}
\exp\left(-H_{N,M}(\mathop{\pmb{\si}})\right) \right) \right].
$$
To this purpose, set
$$
G(\gamma) = \alpha \log \left( \sum_{p=0}^\infty \exp(-\gamma)
\frac{\gamma^p}{p!} \E \left[ \frac{\bar V_{p+1}}{\bar V_p}
\right] \right),$$ where
$$
\bar V_p :=\int \Big \langle \exp\Big( u\Big( \sum_{i \le p} g_{i,M}
\sigma_i \Big)\Big) \Big \rangle_{(x_1,\ldots,x_{p})}
d\mu_{\alpha,\gamma}(x_1) \times \ldots \times
d\mu_{\alpha,\gamma}(x_{p})
$$
and $\langle \cdot \rangle_x$  means integration with respect to
the product measure $\nu$ on $\{-1,1  \}^{p}$
such that $\int \sigma_i d\nu = x_i.$
Then, the main result of this part states that:

\begin{theorem}\label{teopn}
Set $F$ such that $F'(\gamma)=G(\gamma)$ and $F(0)=\log 2 - \alpha
u(0).$ Then, if $\gamma \le \gamma_0$ and (\ref{oblidada}) and
(\ref{e754}) hold true, we have
$$
\vert p_N (\gamma) - F (\gamma) \vert \le \frac{K}{N},
$$
where $K$ does not depend on $\gamma$ and $N$.
\end{theorem}
Since $p_N(0)=\log 2 - \alpha u(0),$ the proof of the theorem is a
consequence of the following proposition.
\begin{proposition}\label{proppn}
If $\gamma \le \gamma_0$  and (\ref{oblidada}) and (\ref{e754})
hold, we have
$$
\vert p'_N (\gamma) - G (\gamma) \vert \le \frac{K}{N},
$$
where $p'_N (\gamma)$ is the right derivative of  $p_N (\gamma)$.
\end{proposition}
\bigskip

\noindent {\bf Proof:} We divide the proof into two steps.

\medskip

\noindent {\it Step 1:} We will check that
\begin{equation}\label{step1}
\vert p'_N(\gamma)- G^1(\gamma) \vert \le \frac{K}{N}.
\end{equation}
where $G^1(\gamma)$ is defined as
$$\alpha \E \left[ \log
\Bigg\langle \exp\left( u \Big( \sum_{i \le N} g_{i,M} \gamma_{i,M}
\sigma_i + g_{N,M} \sigma_N \Big) -  u \Big( \sum_{i \le N} g_{i,M}
\gamma_{i,M} \sigma_i \Big)\right) \Bigg\rangle \right].$$

Following the method used in Lemma 7.4.11 in \cite{Tb}, we
introduce the Hamiltonians
\begin{eqnarray*}
-H_{N,M}^1 (\mathop{\pmb{\si}})= \sum_{k \le M}  u \left( \sum_{i
\le N} g_{i,k}\ (\ga_{i,k} + \delta_{i,k})\ \si_{i} \right),\\
-H_{N,M}^2 (\mathop{\pmb{\si}})= \sum_{k \le M}  u \left( \sum_{i
\le N} g_{i,k}\ \min(1,(\ga_{i,k} + \delta_{i,k}))\ \si_{i}
\right),
\end{eqnarray*}
where $\{ \delta_{i,k} \}_{1\le i \le N, 1 \le k \le M }$ is a
family of i.i.d. random variables with
$P(\delta_{i,k}=1)=\frac{\delta}{N},
P(\delta_{i,k}=0)=1-\frac{\delta}{N}$. We also  assume that this
sequence is independent of all the random sequences previously
introduced. Observe that  the random variables $\min(1,(\ga_{i,k}
+ \delta_{i,k}))$ are
 i.i.d with Bernoulli law of  parameter
 $\frac{\gamma'}{N}$, where
$\gamma'\equiv\gamma+\delta - \frac{\gamma \delta}{N}.$ Set now,
for $j=1,2$,
$$
p_N^j(\delta) = \frac{1}{N} \E \left[ \log \left(
\sum_{\mathop{\pmb{\si}}\in \Sigma_N}
\exp\left(-H_{N,M}^j(\mathop{\pmb{\si}})\right) \right) \right].$$
Obviously, $p_N^2 (\delta)=p_N(\gamma')$, and our first task will
be to show that $p_N^1 (\delta)-p_N^2 (\delta)$ is of order
$\delta^2$: notice that
$$
p_N^1(\delta) - p_N^2(\delta) = \frac{1}{N} \E \left[ \log \langle
\exp\left(-H_{N,M}^1(\mathop{\pmb{\si}}) +
H_{N,M}^2(\mathop{\pmb{\si}}) \right) \rangle_2 \right],$$
where
$\langle \cdot \rangle_2$ denotes the average for the Gibbs' measure
defined by the
Hamiltonian $H_{N,M}^2$. Consider now $Y_{N,M}^1 = \sum_{i,k}
\gamma_{i,k} \delta_{i,k}.$ Since,
$\ga_{i,k}+\delta_{i,k}=\min(1,\ga_{i,k}+\delta_{i,k})+
\ga_{i,k}\delta_{i,k}$, on the set $\{Y_{N,M}^1=0\}$,
we have
$H_{N,M}^1=H_{N,M}^2$. So, we can write
\begin{eqnarray*}
p_N^1(\delta) - p_N^2(\delta) & = & \frac{1}{N} \E \left[
\1_{\{Y_{N,M}^1=1\}} \log \langle
\exp\left(-H_{N,M}^1(\mathop{\pmb{\si}}) +
H_{N,M}^2(\mathop{\pmb{\si}}) \right) \rangle_2 \right] \\
 &  & +  \frac{1}{N} \E \left[
\1_{\{Y_{N,M}^1 \ge 2\}} \log \langle
\exp\left(-H_{N,M}^1(\mathop{\pmb{\si}}) +
H_{N,M}^2(\mathop{\pmb{\si}}) \right) \rangle_2 \right].
\end{eqnarray*}
Using that
$$P(Y_{N,M} \ge 2 ) = 1 - \left( 1 - \frac{\gamma \delta}{N^2} \right)^{NM}
- NM \left( 1 - \frac{\gamma \delta}{N^2} \right)^{NM-1}
\frac{\gamma \delta}{N^2} \le \alpha^2 \delta^2 \gamma^2$$ and
$$ P(Y_{N,M}^1=1)= NM \left(1 - \frac{\gamma \delta}{N^2}
\right)^{NM-1} \frac{\gamma \delta}{N^2} \le \alpha \gamma
\delta,$$
it is easily checked that
\begin{equation}\label{step1001}
\lim_{\delta \longrightarrow 0^+ } \frac{p_N^1(\delta) -
p_N^2(\delta)}{\delta} \le \frac{K}{N},
\end{equation}
which means that we can evaluate the difference
$p_N^1(\delta) - p_N(\gamma)$ instead of
$p_N^2(\delta) - p_N(\gamma)$.

However, following the same arguments as above, we can write
$$
p_N^1(\delta) - p_N(\gamma) = \frac{1}{N} \E \left[ \log \langle
\exp\left(-H_{N,M}^1(\mathop{\pmb{\si}}) +
H_{N,M}(\mathop{\pmb{\si}}) \right) \rangle \right].$$ We consider
now $Y_{N,M} = \sum_{i,k} \delta_{i,k}.$ Notice that on the set
$\{Y_{N,M}=0\}$, $H_{N,M}=H_{N,M}^1$. So, we can write
\begin{eqnarray*}
p_N^1(\delta) - p_N(\gamma) & = & \frac{1}{N} \E \left[
\1_{\{Y_{N,M}=1\}} \log \langle
\exp\left(-H_{N,M}^1(\mathop{\pmb{\si}}) +
H_{N,M}(\mathop{\pmb{\si}}) \right) \rangle \right] \\
 &  & +  \frac{1}{N} \E \left[
\1_{\{Y_{N,M} \ge 2\}} \log \langle
\exp\left(-H_{N,M}^1(\mathop{\pmb{\si}}) +
H_{N,M}(\mathop{\pmb{\si}}) \right) \rangle \right] \\
&\equiv& V_1(\delta) + V_2(\delta).
\end{eqnarray*}

Let  us bound now $V_1(\delta)$ and $V_2(\delta)$: since
$$P(Y_{N,M} \ge 2 ) = 1 - \left( 1 - \frac{\delta}{N} \right)^{NM}
- NM \left( 1 - \frac{\delta}{N} \right)^{NM-1} \frac{\delta}{N}
\le \alpha (NM-1) \delta^2,$$ we have
$$
\vert V_2 (\delta) \vert \le 2 \alpha^2 (NM-1) U_\infty
\delta^2.$$ On the other hand, using a symmetry argument,
we get
\begin{align*}
&  V_1 (\delta) =  NM \left(1 - \frac{\delta}{N} \right)^{NM-1}
\frac{\delta}{N^2} \\
&  \times \E \left[ \log \Bigg\langle
\exp\left( u \Big( \sum_{i \le N} g_{i,M} \gamma_{i,M} \sigma_i +
g_{N,M} \sigma_N \Big) -  u \Big( \sum_{i \le N} g_{i,M}
\gamma_{i,M} \sigma_i \Big)\right)
\Bigg\rangle \right].
\end{align*}
Hence, we obtain that
\begin{align}
& \lim_{\delta \longrightarrow 0^+ } \frac{ p_N^1(\delta) -
p_N(\gamma)}{\delta} = \lim_{\delta \longrightarrow 0^+ } \frac{
V_1(\delta) +
V_2(\delta)}{\delta}\label{step1002}\\
&= \alpha \E \left[ \log \Bigg\langle \exp\left( u \Big( \sum_{i
\le N} g_{i,M} \gamma_{i,M} \sigma_i + g_{N,M} \sigma_N \Big) -  u
\Big( \sum_{i \le N} g_{i,M} \gamma_{i,M} \sigma_i \Big)\right)
\Bigg\rangle \right].\nonumber
\end{align}
Eventually, since
$$
p_N'(\gamma)= \lim_{\gamma' \longrightarrow \gamma^+}
\frac{p_N(\gamma')-p_N(\gamma)}{\gamma'-\gamma} = \lim_{\delta
\longrightarrow 0^+} \frac{p_N^2 (\delta)-p_N(\gamma)}{\delta
\left( 1 - \frac{\gamma}{N} \right)},$$ putting together
(\ref{step1001}) and (\ref{step1002}), we obtain (\ref{step1}).

\medskip

\noindent {\it Step 2:} Let us check now that
\begin{equation}\label{step2}
\vert G(\gamma)- G^1(\gamma) \vert \le \frac{K}{N}.
\end{equation}
To this purpose, set
$$
\Psi :=  \Bigg\langle \exp\left( u \Big( \sum_{i \le N} g_{i,M}
\gamma_{i,M} \sigma_i + g_{N,M} \sigma_N \Big) -  u \Big( \sum_{i \le N}
g_{i,M} \gamma_{i,M} \sigma_i \Big)\right) \Bigg\rangle ,
$$
and let us try to evaluate first $\E[\Psi]$: notice  that
\begin{eqnarray*}
\Psi
& = & \frac{\sum_{\mathop{\pmb{\si}}\in \Sigma_N}
 \exp\left( u ( \sum_{i \le N} g_{i,M}
\gamma_{i,M} \sigma_i + g_{N,M} \sigma_N ) \right)
\exp\left(-H_{N,M-1}(\mathop{\pmb{\si}})\right)}{\sum_{\mathop{\pmb{\si}}\in
\Sigma_N} \exp\left(-H_{N,M}(\mathop{\pmb{\si}})\right)} \\
& = & \frac{\langle
 \exp\left( u ( \sum_{i \le N} g_{i,M}
\gamma_{i,M} \sigma_i + g_{N,M} \sigma_N ) \right)\rangle_{M-1}}
{\langle \exp\left( u ( \sum_{i \le N} g_{i,M} \gamma_{i,M}
\sigma_i) \right)\rangle_{M-1}},
\end{eqnarray*}
where $\langle \cdot \rangle_{M-1}$ denotes the usual average
using the Hamiltonian $H_{N,M-1}$. Set $B_p:=\{\sum_{i=1}^{N-1}
\gamma_{i,M}=p, \gamma_{N,M}=0\}$ and $B:=\{\gamma_{N,M}=1\}$, and
let us denote by $\E_M$ the conditional expectation given
$\{\gamma_{i,M}, 1 \le i \le N\}$. Then
\begin{eqnarray}
& &\E \left[ \Psi \right] = \E \left[ \sum_{p=0}^{N-1} \1_{B_p}
\E_M \left[ \Psi  \right] \right] + \E \left[ \1_{B} \E_M \left[
\Psi
 \right] \right] \label{step2000}\\
& & = \sum_{p=0}^{N-1} \binom{N-1}{p} \left( \frac{\gamma}{N}
\right)^p \left( 1- \frac{\gamma}{N} \right)^{N-p+1} \E \left[
\frac{\langle \exp(V_{p+1}) \rangle_{M-1}}{\langle \exp(V_{p})
\rangle_{M-1}} \right] + \frac{\gamma}{N} e^{2 U_\infty},\nonumber
\end{eqnarray}
where $$V_p:=u \left(\sum_{i \le p} g_{i,M} \sigma_i \right).$$
Set ${\mathop{\pmb{X}_{p}}}=(\langle \si_1 \rangle, \ldots,
\langle \si_p \rangle)$. Then, using the triangular inequality and
following the same arguments as in Proposition \ref{p727}, we get,
for a strictly positive constant $K$,
\begin{eqnarray}
& &\left \vert \E \left[ \frac{\langle \exp(V_{p+1})
\rangle_{M-1}}{\langle \exp(V_{p}) \rangle_{M-1}} \right] - \E
\left[ \frac{\langle \exp(V_{p+1})
\rangle_{{\mathop{\pmb{X}_{p+1}}}}}{\langle
\exp(V_{p}) \rangle_{{\mathop{\pmb{X}_{p}}}}} \right] \right\vert \nonumber\\
& & = \left \vert \E \left[ \frac{\langle \exp(V_{p+1})
\rangle_{M-1} \langle \exp(V_{p}) \rangle_{{\mathop{\pmb{X}_{p}}}}
- \langle \exp(V_{p+1}) \rangle_{{\mathop{\pmb{X}_{p+1}}}} \langle
\exp(V_{p}) \rangle_{M-1} }{\langle \exp(V_{p}) \rangle_{M-1}
\langle \exp(V_{p})
\rangle_{{\mathop{\pmb{X}_{p}}}}} \right]\right\vert \nonumber \\
& &\le e^{3 U_\infty} \E \left[ \left \vert \langle \exp(V_{p+1})
\rangle_{M-1} - \langle \exp(V_{p+1})
\rangle_{{\mathop{\pmb{X}_{p+1}}}} \right \vert \right. \nonumber \\
& & \qquad \left. + \left \vert \langle \exp(V_{p}) \rangle_{M-1}
- \langle \exp(V_{p}) \rangle_{{\mathop{\pmb{X}_{p}}}} \right
\vert\right]
\nonumber \\
& & \le e^{3 U_\infty} \frac{p^2 K}{N}. \label{step2001}
\end{eqnarray}
Consider now some i.i.d. random variables $z_1, \ldots, z_p$ of
law $\mu_{\alpha,\gamma}$ such that (\ref{limmag}) holds. Set
${\mathop{\pmb{Y}_{p}}}=(z_1,\ldots,z_p).$ Then, following the
same arguments as above, we get, for a strictly positive constant
$K$,
\begin{eqnarray}
& &\left \vert \E \left[ \frac{\langle \exp(V_{p+1})
\rangle_{{\mathop{\pmb{X}_{p+1}}}}}{\langle \exp(V_{p})
\rangle_{{\mathop{\pmb{X}_{p}}}}} \right] - \E \left[
\frac{\langle \exp(V_{p+1})
\rangle_{{\mathop{\pmb{Y}_{p+1}}}}}{\langle
\exp(V_{p}) \rangle_{{\mathop{\pmb{Y}_{p}}}}} \right] \right\vert \nonumber\\
& &\le e^{3 U_\infty} \E \left[ \left \vert \langle \exp(V_{p+1})
\rangle_{{\mathop{\pmb{X}_{p+1}}}} - \langle \exp(V_{p+1})
\rangle_{{\mathop{\pmb{Y}_{p+1}}}} \right \vert \right. \nonumber
\\ & & \qquad \left. + \left \vert \langle \exp(V_{p})
\rangle_{{\mathop{\pmb{X}_{p}}}} - \langle
\exp(V_{p}) \rangle_{{\mathop{\pmb{Y}_{p}}}} \right \vert\right] \nonumber\\
& & \le e^{3 U_\infty}\frac{p^3 K}{N},\label{step2002}
\end{eqnarray}
where in the last inequality we have used (\ref{limmag})  and
the fact that
$$
\frac{\partial}{\partial x_i} \langle \exp(V_p) \rangle_x \le
e^{U_\infty}.
$$

Notice that if $W$ is a random variable with law
Bin$(N-1,\frac{\gamma}{N})$, then $\E(W^3) \le K$, where $K$ does
not depend on $N$. So, putting together (\ref{step2000}),
(\ref{step2001}) and (\ref{step2002}), we get
$$
\E \left[ \Psi \right] =\sum_{p=0}^{N-1}  \binom{N-1}{p} \left(
\frac{\gamma}{N} \right)^p \left( 1- \frac{\gamma}{N}
\right)^{N-p} \E \left[ \frac{\langle \exp(V_{p+1})
\rangle_{\mathop{\pmb{Y}_{p+1}}}}{\langle \exp(V_{p})
\rangle_{\mathop{\pmb{Y}_{p}}}} \right] + \frac{K}{N}.
$$
Using now similar arguments to those ones used in the proof of
Lemma \ref{l745}, we get
\begin{equation}\label{hk}
\E \left[ \Psi \right] =\sum_{p=0}^{\infty} \exp(-\gamma)
\frac{\gamma^p}{p!} \E \left[ \frac{\langle \exp(V_{p+1})
\rangle_{{\mathop{\pmb{Y}_{p+1}}}}}{\langle \exp(V_{p})
\rangle_{{\mathop{\pmb{Y}_{p}}}}} \right] + \frac{K}{N}.
\end{equation}
Eventually, once (\ref{hk}) has been obtained,
(\ref{step2}) can be established following the method used
in Proposition 7.4.10 in \cite{Tb}, the remaining details being
left to the reader.
 \hfill$\square$

\end{document}